\documentclass[a4paper,11pt,oneside]{article}

\usepackage{amsmath,amsthm,amsfonts,amssymb,amscd}
\usepackage{subfigure}
\usepackage{booktabs}
\usepackage{graphicx,xcolor}
\usepackage{algorithm}
\usepackage{algpseudocode}
\usepackage[ngerman,english]{babel}
\usepackage{siunitx}

\newcommand\plot[1]{\let\frame\relax% Use \fbox to show frame.
  \frame{\includegraphics[clip,trim=0 220 0 60,width=8cm]{#1}}}
\newcommand\high{^{\textup{high}}}
\newcommand\tsb[1]{_{\textup{#1}}}
\newcommand\Ieff{I\tsb{eff}}
\newcommand\Iind{I\tsb{ind}}
\newcommand\Jdrag{J\tsb{drag}}

\usepackage{authblk}
\usepackage[top=2cm, bottom=2cm, left=2cm, right=2cm]{geometry}
\theoremstyle{plain} % default
\newtheorem{theorem}{Theorem}[section]
\newtheorem{Proposition}{Proposition}[section]

\newtheorem{Corollary}[theorem]{Corollary}

\newtheorem{Definition}[theorem]{Definition}

\newtheorem{Problem}[theorem]{Problem}

\theoremstyle{definition} %

\theoremstyle{remark} %

\usepackage{xcolor}

\newcommand{\bs}[1]{{#1}} % use normal

\let\Hat\widehat
\let\widehat\hat

\newcommand{\alemap}{\hat{\cal A}}
\newcommand{\hOmega}{\widehat\Omega}

\newcommand{\hGamma}{\widehat\Gamma}

\newcommand{\hX}{\bs{\widehat{X}}}

\newcommand{\hV}{\bs{\widehat{V}}}

\newcommand{\hF}{\bs{\widehat{F}}}
\newcommand{\hJ}{\widehat{J}}
\newcommand{\hE}{\bs{\widehat{E}}}

\newcommand{\hdiv}{\Hat{\operatorname{div}}}

\newcommand{\tr}{{\textup{tr}}}

\newcommand{\hnabla}{\bs{\widehat\nabla}}
\newcommand{\hSigma}{\bs{\widehat\Sigma}}

\newcommand{\hsigma}{\bs{\widehat\sigma}}

\newcommand{\hrho}{\bs{\widehat\rho}}

\newcommand{\hv}{\bs{\widehat v}}
\newcommand{\hu}{\bs{\widehat u}}
\newcommand{\hp}{\bs{\widehat p}}

\newcommand{\hf}{\bs{\widehat f}}
\newcommand{\hg}{\bs{\widehat g}}

\newcommand{\hA}{\bs{\widehat A}}

\newcommand{\hU}{\bs{\widehat U}}

\newcommand{\hz}{\bs{\widehat z}}
\newcommand{\hZ}{\bs{\widehat Z}}

\newcommand{\hPsi}{\bs{\widehat\Psi}}
\newcommand{\hPhi}{\bs{\widehat\Phi}}
\newcommand{\hpsi}{\bs{\widehat\psi}}

\newcommand{\hn}{\bs{\widehat n}}

\newif\ifMAKEPICS
\MAKEPICStrue

\ifMAKEPICS
\usepackage[cleanup,subfolder]{gnuplottex}
\usepackage{xparse}

\ExplSyntaxOn
\DeclareExpandableDocumentCommand{\convertlen}{ O{cm} m }
{
	\dim_to_decimal_in_unit:nn { #2 } { 1 #1 } cm
}
\ExplSyntaxOff
\fi

\begin{document}

\title{
%Multiple goal-oriented error estimation
%for the stationary p-Navier-Stokes equations
% 2nd choice for the future
Multigoal-oriented error estimation and mesh adaptivity
for fluid-structure interaction
}

\author[1]{K. Ahuja}
\author[2,3,5,6]{B. Endtmayer}
\author[2]{M. C. Steinbach}
\author[2,3,4]{T. Wick}

%\author[2]{Bernhard Endtmayer}
%\author[1]{Marc Steinbach}
%\author[3]{Kapil Ahuja}
%\author[1]{Thomas Wick}

\affil[1]{Discipline of Computer Science and Engineering,
Indian Institute of Technology Indore, Indore, India}

\affil[2]{Leibniz Universit\"at Hannover, Institut f\"ur Angewandte
  Mathematik, Welfengarten 1, 30167 Hannover, Germany}

\affil[3]{Cluster of Excellence PhoenixD (Photonics, Optics, and
	Engineering -- Innovation Across Disciplines), Leibniz Universit\"at Hannover, Germany}

%\affil[2]{Johann Radon Institute for Computational and Applied Mathematics,
%Austrian Academy of Sciences, Altenbergerstr. 69, A-4040 Linz, Austria}

\affil[4]{Universit\'e Paris-Saclay, ENS Paris-Saclay, LMT -- Laboratoire de M\'ecanique et Technologie, 91190 Gif-sur-Yvette, France}

\affil[5]{ Johannes Kepler University Linz,
		Doctoral Program on Computational Mathematics, Altenbergerstr. 69, A-4040 Linz, Austria}
\affil[6]{Austrian Academy of Sciences, Johann Radon Institute for Computational and Applied Mathematics,  Altenbergerstr. 69, A-4040 Linz, Austria}

\date{}

\maketitle

%%%%%%%%%%%%%%%%%%%%%%%%%%%%%%%%%%%%%%%%%%%%%%%%%%%%%%%%%%%%
\begin{abstract}
In this work, we consider multigoal-oriented
error estimation for stationary fluid-structure interaction.
The problem is formulated within a variational-monolithic
setting using arbitrary Lagrangian-Eulerian coordinates.
Employing the dual-weighted residual method for goal-oriented
a posteriori error estimation, adjoint sensitivities are required.
For multigoal-oriented error estimation, a combined functional
is formulated such that
several quantities of interest are controlled simultaneously.
As localization technique for mesh refinement we employ a partition-of-unity.
Our algorithmic developments are substantiated with several
numerical tests such as an elastic lid-driven cavity with two goal functionals,
an elastic bar in a chamber with two goal functionals, and the FSI-1 benchmark
with three goal functionals.
\end{abstract}

%%%%%%%%%%%%%%%%%%%%%%%%%%%%%%%%%%%%%%%%%%%%%%%%%%%%%%%%%%%%
\section{Introduction}
\label{sec_intro}
This work is devoted to
multigoal-oriented a posteriori error estimation
and mesh adaptivity for the stationary fluid-structure
interaction. Therein, the incompressible Navier-Stokes equations
interact with elastic solids. References to extensive work
on fluid-structure interaction include the textbooks \cite{BuSc06,FoQuaVe09,GaRa10,BuSc10,BaTaTe13,BoGaNe14,Ri17_fsi,FrHoRiWiYa17}.

Multigoal-oriented a posteriori error estimation
using the dual-weighted residual method \cite{BeRa01,BaRa03} with adjoint sensitivity measures
goes back to Hartmann/Houston \cite{HaHou03} and Hartmann
\cite{Ha08}. The principal idea is to construct
a combined goal functional $J_c$
as a linear combination of the individual goal functionals.
The adjoint problem is then solved with respect to
this combined goal functional.

In recent years, the multigoal-oriented error estimates have become
very attractive as attested in several studies
\cite{PARDO20101953,BruZhuZwie16,KerPruChaLaf2017,EndtNeiLaWiWo20}.
Moreover, several improvements in the algorithmic
techniques have been undertaken.
In \cite{EndtWi17}, the sign computation was facilitated
such that only one additional problem of higher-order must be solved.
Then, in \cite{EndtLaWi18}, the error estimator was extended
to balance the nonlinear iteration and the discretization
error. For single goal functional this approach
was first developed in \cite{RanVi2013}.
Recently, two-sided error estimates could be proven
under a saturation assumption \cite{EndtLaWi20}; a summary
is contained in the PhD thesis of the second author \cite{Endt21}.
Proofs for the saturation assumption are given in
\cite{DoerflerNochetto2002, Agouzal2002}. For problems containing convection
or transport, the assumption is not always true \cite[Example 7.2]{EndtLaWi20}.
%For elliptic problems
%this assumption is satisfied, while for problems containing
%convection or transport, the assumption is not always true
%\cite[Example 7.2]{EndtLaWi20}.
This
has consequences in this paper, since we consider again
a system with convection terms in the Navier-Stokes part. Therefore,
we cannot expect asymptotically perfect effectivity indices,
but nonetheless the algorithms are reliable.

The objective in this paper is to apply the multigoal-oriented
framework developed in \cite{EndtLaWi18} to fluid-structure interaction, which
is a multidomain multiphysics problem.
This is an illustrative example of a multiphysics problem in which various quantities
in the different subdomains (flow and solid) might be controlled simultaneously. Specific interest
is on drag, lift, pressure and displacement evaluations.
Our model is based on the ALE (arbitrary Lagrangian-Eulerian) framework
\cite{DoFaGi77,HuLiZi81,FoNo99} using variational-monolithic coupling. Based
on this monolithic model and combined goal functional, the Lagrangian is defined
from which the primal model and the adjoint equation are derived. The resulting
a posteriori error estimator is then localized using a partition-of-unity
technique \cite{RiWi15_dwr}, which is specifically attractive for
multidomain multiphysics problems, because the error estimator is directly
based on the weak (variational) formulation without the traditional backward
partial integration yielding that would yield second-order operators and interface terms to be
evaluated. Consequently, the partition-of-unity is less error-prone
as well as more efficient. Based on the localized error estimator, an adaptive
scheme is derived from which error-controlled mesh adaptivity is carried out
and for which two well-known numerical tests, the FSI-1 benchmark
\cite{HrTu06b,BuSc06,BuSc10} and a modified lid-driven cavity example with
elastic volume at the bottom \cite{Du06,Du07,RiWi10}
are taken to substantiate our algorithmic developments.

The outline of this paper is as follows. In Section
\ref{sec_problems} we formulate the primal
problem and its Galerkin finite element discretization.
In Section \ref{sec_MG}, multigoal-oriented error estimates
are derived. Then, in Section \ref{sec_error_loc}
the errors are localized to carry out adaptive mesh refinement.
In Section \ref{sec_tests}, three
numerical tests are investigated in order to show
the possibilities using our proposed multigoal-oriented
framework. In the conclusions, our work is summarized.

%%%%%%%%%%%%%%%%%%%%%%%%%%%%%%%%%%%%%%%%%%%%%%%%%%%%%%%%%%%%
\section{Fluid-structure interaction: primal problem and discretization}
\label{sec_problems}
In this section, we introduce the governing equations and
formulate the adjoint problem. The latter is required
for goal-oriented a posteriori error estimation with
the dual-weighted residual (DWR) method \cite{BeRa01}.

\subsection{Notation}
We denote by $\Omega \subset \mathbb{R}^d$, $d=2$, the domain of the
fluid-structure interaction problem. This domain consists of the subdomains
$\Omega_f$ and $\Omega_s$. The interface between both
domains is denoted by $\Gamma_i  = \partial\Omega_f  \cap \partial\Omega_s$.
The reference
domains (obtained from a change of the coordinate system)
are denoted by $\hOmega_f$ and $\hOmega_s$, respectively, with the
interface $\hGamma_i$.
Further, we denote the outer boundary with
$\partial\hOmega = \hGamma = \hGamma_D \cup \hGamma_N$
where
$\hGamma_D $ and $\hGamma_N $ denote Dirichlet and Neumann boundaries, respectively.

We adopt standard notation for the usual Lebesgue and Sobolev spaces \cite{Wlo87}.
We use the notation $(\cdot ,\cdot)$ for a scalar product
on a Hilbert space $X$ on the domain $\Omega$
and $\langle \cdot, \cdot \rangle$ for the scalar product on the boundary
$\partial \Omega$. Finally, we introduce
short hand notations for the Sobolev spaces:
%For the time independent functions in our static problem
%the Sobolev spaces are defined as
$L := L^2 (\Omega)$ and
$V := H^1 (\Omega)$, respectively
$\hat L := L^2 (\hat\Omega)$ and
$\hat V := H^1 (\hat\Omega)$. If the fluid and solid domains
must be distinguished, we add subscripts $f$ or $s$.
For Dirichlet boundaries $\hGamma_D $ with homogeneous
conditions, we use $V^0$ and $\hat V^0$, respectively.

\subsection{Fluid flow in ALE coordinates}
Let $\alemap_f\colon \hOmega_f \to \Omega_f$
be an invertible mapping, the so-called
ALE transformation \cite{DoFaGi77,HuLiZi81,FoNo99}.
We define the unknowns
$\hat v_f$ and $\hat p_f$ in $\hat\Omega_f$ by
\begin{equation*}
\hat v_f(\hat x) = v_f(x) = v_f(\alemap_f(\hat x)),\quad
\hat p_f(\hat x) = p_f(x) = p_f(\alemap_f(\hat x)).
\end{equation*}
Then, with
\begin{equation*}
\hat F_f:=\hat\nabla\alemap_f,\quad \hat J_f:=\det(\hat F_f),
\end{equation*}
we obtain the relations
\begin{equation*}
\nabla v_f = \hat\nabla\hat v_f\hat F_f^{-1},\quad
%\partial_t v_f = \partial_t \hat v_f - (\hat F_f^{-1} \partial_t
%\alemap_f\cdot\hat\nabla )\hat v_f,\quad
\int\limits_{\Omega_f} f(x)\,\text{d}x = \int\limits_{\hat\Omega_f} \hat
f(\hat x) \hat J_f\,\text{d}\hat x.
\end{equation*}
With these relations, we can formulate the Navier-Stokes equations in artificial
coordinates:
\begin{Problem}[Variational fluid problem, ALE framework]
  \label{eq:fluid:ale}
Let $\hat v_f^D$ denote non-homogeneous flow Dirichlet boundary conditions.
Find $\{\hat v_f,\hat p_f\} \in \{ \hat v_f^D + \hat V_f\} \times
\hat {L}_f$, such that
  \begin{equation*}
    \begin{aligned}
      (\hat J_f \hat\rho_f  (\hat F_f^{-1}\hat
      v_f\cdot\hat\nabla) \hat
      v_f,\hat\psi^v)_{\hat\Omega_f}
      + (\hat J_f\hat\sigma_f \hat
      F_f^{-T},\hat\nabla\hat\psi^v)_{\hat\Omega_f}
      - \langle \hat g_f
      ,\hat \psi^v\rangle_{\hat\Gamma_i\cup\hat\Gamma_N}
      & = 0
      &&\forall\hat\psi^v\in\hat V_f^0,
      \\
      (\hdiv\,(\hat J_f\hat F_f^{-1}\hat
      v_f),\hat\psi^p)_{\hat\Omega_f}
      &=0&&\forall\hat\psi^p\in\hat L_f^0,
    \end{aligned}
  \end{equation*}
  with the ALE transformed Cauchy stress tensor
  \begin{equation*}
  \hat \sigma_f := -\hat p_fI +\hat\rho_f\nu_f(\hat\nabla \hat v_f \hat F_f^{-1}
  + \hat F_f^{-T}\hat\nabla \hat v_f^T).
  \end{equation*}
  The viscosity and the density of the fluid are denoted by $\nu_f$
and $\hat\rho_f$, respectively. The term $I$ denotes the
identity matrix in $\mathbb{R}^{d\times d}$.
The function $\hat g_f$ represents
Neumann boundary conditions for both physical boundaries (e.g.,
stress zero at outflow boundary),
and normal stresses on $\hat\Gamma_i$.
We note that the specific choice of the transformation $\alemap_f$ is
up to now arbitrary and left open.
\end{Problem}

\subsection{Solid description in Lagrangian coordinates}
Solids are usually formulated in Lagrangian
coordinates, which means that one needs to find a mapping from the physical (deformed) domain
$\Omega_s$ to the reference domain $\hat\Omega_s$.

The transformation $\alemap_s\colon\hOmega_s \to\Omega_s$
is naturally given by the deformation itself:
\vskip-.6cm
\begin{eqnarray}
\alemap_s(\hat x) = \hat x+\hat u_s(\hat x),\quad
\hat F_s:=\hat\nabla\alemap_s = I +\hat\nabla\hat u_s,\quad
\hat J_s:=\det(\hat F_s ).
\end{eqnarray}
This transformation is identical in its definition
to the ALE transformation $\alemap_f$, and for this reason
both can be identified by omitting the subscripts $f$ and $s$.

As material law, we work with
the elastic compressible
(geometrically) nonlinear Saint Venant-Kirchhoff material (STVK).
It is well suited for (relatively) large displacements with the limitation
of small strains.
The strain is defined by $\hat E:=\frac{1}{2}(\hat F^T \hat F -I)$.
The sought physical unknowns are the displacement $\hu$ and the
velocity $\hv$:

\begin{Problem}[Solid in Lagrangian coordinates]
  \label{eq:stvk:lagrange}
  Find $\hat u_s \in \hat {V}^0_s$
  such that
  \begin{eqnarray}
    \begin{aligned}
      %(\hat\rho_s \partial_t\hat v_s,\hat\psi^v)_{\hat\Omega_s}
      (\hat J_s \hat\sigma_s \hat F_s^{-T},\hat\nabla\hat\psi^v)_{\hat\Omega_s}
      -\langle\hat J_s \hat\sigma_s \hat
      F_s^{-T} \hat n_s ,\hat\psi^v  \rangle_{\hat\Gamma_i\cup\hat\Gamma_N}
      &= (\hat\rho_s \hat f_s , \hat \psi^v )_{\hat\Omega_s}
      &&\quad\forall\hat \psi^v\in\hat V_s^0,\\
      %%%%%%%%%%%%%%%%
      %(\partial_t\hat u_s-\hat v_s,\hat\psi^u)_{\hat\Omega_s}&=0
      %&&\quad\forall\hat \psi^u\in\hat V_s    ,\\
    \end{aligned}
  \end{eqnarray}
  where $\hat n_s$ is the outer normal vector
  on $\hat\Gamma_i$ and $\hat\Gamma_N$, respectively.
  The Cauchy stress tensors for STVK material is
given by
  \begin{align}
    \hat\sigma_s &:=\hat J^{-1} \hat F (\lambda_s (\tr\hat E) I
    + 2\mu_s \hat E)\hat F^{T} \label{STVK_material}
  \end{align}
with the positive Lam\'e coefficients
$\mu_s$ and $\lambda_s$. External volume forces
are described by the term $\hat f_s$.
\end{Problem}

%%%%%%%%%%%%%%%%%%%%%%%%%%%%%%%%%%%%%%%%%%%%%%%%%%

\subsection{The coupled problem in ALE coordinates}\label{sec:eq:fsi}
In the solid part, $\alemap_s$ is determined naturally
by the solid displacements (see for example \cite{Cia1984}).
Therefore, only the flow part $\alemap_f$ needs to be specified.
On the interface $\hGamma_i$,
this transformation is given by taking the solid displacement
$\hat u_f = \hat u_s$ such that we can define
\begin{eqnarray}
\alemap_f(\hat x)\big|_{\hat\Gamma_i} :=
\hat x+\hat u_s(\hat x)\big|_{\hat\Gamma_i}.
\end{eqnarray}
On the outer boundary of the fluid domain,
$\partial\hOmega_f\setminus\hGamma_i$, there holds
$\alemap_f=\text{id}$. Inside $\hOmega_f$, the transformation should
be sufficiently smooth and regular.
For some relevant mathematical regularity results,
we refer to \cite{FoNo99,FoNo04}.
Away from the interface $\hat\Gamma_i$ the mapping
can be extended arbitrarily in various fashions, for example
harmonic, nonlinear harmonic/elastic \cite{StTeBe03},
or biharmonic \cite{Wi11}, to mention just a few cases.
Some recent results on mesh motion can be found in \cite{MR4222508}.
The simplest harmonic model reads (in strong formulation):
\begin{Problem}[Harmonic mesh motion]
\label{problem_harmonic}
Find $ \hat u_f\colon\hOmega_f\to\mathbb{R}^d$ such that
\begin{equation}\label{eq:cauchy_stress_harmonic}
\hdiv (\hat\sigma_g ) = 0,\quad
\hat u_f=\hat u_s\text{ on }\hat\Gamma_i,\quad
\hat u_f=0\text{ on }\partial\hat\Omega_f\setminus\hat\Gamma_i,
\end{equation}
with $\hat\sigma_g = \hat\alpha_u \hat\nabla \hat u_f $ and
$\hat\alpha_u > 0$.
\end{Problem}

In the following,
we focus on a variational-monolithic description of the coupled problem
\cite{Ri17_fsi,Wi20_book}.
To this end, we define a continuous variable $\hat u$ in $\hat\Omega$
defining the deformation in $\hat\Omega_s$ and supporting the
transformation in $\hat\Omega_f$. Thus, we drop subscripts on $\hu$, and
because the definition of $\alemap_f$ coincides with the previous
definition of $\alemap_s$, we define in $\hat\Omega$:
\begin{eqnarray}
\alemap:=\text{id}+\hat u,\quad
\hat F:=I+\hat\nabla \hat u,\quad
\hat J:=\det(\hat F).
\end{eqnarray}
Then, we obtain the overall problem formulation:
\begin{Problem}[Variational-monolithic fluid-structure interaction framework]
  \label{eq:fsi:ale:harmonic}
  $ $\\Find $\{\hat v,\hat u,\hat p\} \in \{ \hat v^D + \hat {V}\}
\times \hat {V} \times \hat {L}$ such that
  \begin{eqnarray*}
    \begin{aligned}
      %(\hat J\hat\rho_f  \partial_t \hat v,\hat\psi^v)_{\hat\Omega_f}
      (\hat\rho_f \hat J  (\hat F^{-1}\hat
      v%-\partial_t \hat u
\cdot\hat\nabla) \hat v),
      \hat\psi^v)_{\hat\Omega_f}
      + (\hat J\hat\sigma_f\hat
      F^{-T},\hat\nabla\hat\psi^v)_{\hat\Omega_f}
      %
      %+ (\hat\rho_s \partial_t \hat v,\hat\psi^v)_{\Omega_s}
      + (\hat J\hat\sigma_s\hat F^{-T},\hat\nabla\hat\psi^v)_{\hat\Omega_s}&\\
      {} - \langle \hat g, \hat\psi^v \rangle_{\hat\Gamma_N} -
      (\hat\rho_f \hat J\hat f_f, \hat\psi^v)_{\hat\Omega_f}
      - (\hat\rho_s\hat f_s, \hat\psi^v)_{\hat\Omega_s}
      &=0&&\forall\hat\psi^v\in \hat {V}^0,
      \\
      %%%%%%%%%%%%
      (%\partial_t\hat u-
      \hat v,\hat\psi^u)_{\hat\Omega_s}
      +
      (\hat\sigma_g ,\hat \nabla\hat\psi^u)_{\hat\Omega_f}
      %-\langle \hat\sigma_g \hat n_f ,\hat\psi^u\rangle_{\hat \Gamma_i}
      &=0&&\forall\hat\psi^u\in \hat {V}^0,\\
      %%%%%%%%%%%%%%%%%%%%%%%%%%%%%%
      (\hdiv\,(\hat J\hat F^{-1}
      \hat v_f),\hat\psi^p)_{\hat\Omega_f}
      %+ (\hat p_s ,\hat \psi^p)_{\hat\Omega_s}
      &=0&&\forall\hat\psi^p\in \hat {L},
    \end{aligned}
  \end{eqnarray*}
  with $\hat\rho_f$, $\nu_f$, $\mu_s$, $\lambda_s$, $\hat F$, and $\hat J$.
  The stress tensors
  $\hat\sigma_f$, $\hat\sigma_s$, and $\hat\sigma_g$ are defined in
  Problems \ref{eq:fluid:ale}, \ref{eq:stvk:lagrange},
  and \ref{problem_harmonic}, respectively.
Written in an abstract form, we have:
find $\hU\in\hX:=\{ \hat v^D + \hat {V}\}
\times \hat {V} \times \hat {L}$ such that
\[
\hA(\hU)(\hPsi) = 0 \quad\forall\hPsi\in\hX^0,
\]
where $\hat X^0 := \hat V^0 \times \hat V^0 \hat \times \hat L$, and
where the semi-linear form $\hA(\hU)(\hPsi)$ is nonlinear in $\hU$ and linear
in the test function $\hPsi$.
\end{Problem}

\subsection{Galerkin finite element discretization}
\label{subsec_SpatialDiscretization}
We assume that $\Omega \subset \mathbb{R}^d$ is a polyhedral domain.
Let $\mathcal{T}_h$ be a subdivision of $\Omega$ into quadrilateral elements
such that $\bigcup_{K \in \mathcal{T}_h}\overline{K}=\overline{\Omega}$
and $K \cap K'= \emptyset$ for all $K,K' \in \mathcal{T}_h$ with $K \neq K'$.
The discretization parameter is denoted by $h$.

Working with adaptive meshes, to ensure global continuity and therefore
global conformity, the degrees of
freedom on interfaces between different elements of
different refinement levels have to fulfill additional constraints.
These are obtained by interpolation where hanging nodes \cite{CaOd84}
do not carry any degrees of freedom.
Specifically, we use continuous
tensor-product finite elements as described in
\cite{Ciarlet:2002:FEM:581834} and  \cite{Braess}.
%We also refer the reader to \cite{ArnBofFal2002}.

Using such a conforming method, we have $\hat X_h \subset \hat X$,
where $\hat X_h$ denotes the finite-dimensional space.
This space is composed of the subspaces $\hat V_h$ and $\hat W_h$
for the unknowns $\hat v_h$ and $\hat p_h$, respectively.
To satisfy discrete inf-sup stability (also known as LBB condition),
we must choose the discrete function spaces in a proper way.
To this end, we work with the well-known Taylor-Hood element
$Q_c^2/Q_c^1$ using biquadratic functions for the velocities
and bilinear functions for the pressure.

The resulting discrete primal problem reads as follows in
a compact abstract notation:
find $\hat U_h \in \hat X_h$ such that
\begin{equation} \label{problem: discrete primal problem}
{\hat A}(\hat U_h)(\hat\Psi_h)=0 \quad \forall \hat\Psi_h \in \hat X_h^0.
\end{equation}

\section{Multigoal-oriented a posteriori error estimation}
\label{sec_MG}
In this section, we describe multigoal-oriented
error estimation with the dual-weighted residual method.
For examples
of goal-oriented error estimation in fluid-structure interaction
we refer the reader to
\cite{VANDERZEE20112738,FickBrummelenZee2010,Ri12_dwr,GraeBa06,Du06,DuRaRi09,Ri17_fsi,FaiWi18}
and \cite{Wi21_git} including an open-source github version for single goal functionals.

\subsection{Optimization problem}

Goal-oriented error estimation aims at controlling the discretization error
between the continuous solution $\hat U\in \hat X$
and the discrete solution $\hat U_h\in \hat X_h$
measured in terms of a goal functional $J\colon \hat X\to\mathbb{R}$.
Such goal functionals describe technical quantities of interest such as
drag or lift values, or point evaluations, but also include global norm errors.

In order to derive sensitivity information, the evaluation of $J(\hat U)$
is artificially reformulated as an optimization problem
where the constraint already determines $\hat U$ uniquely \cite{BeRa01},
\[
  \min \, J(\hat U) \quad\text{s.t.}\quad
  \hat A(\hat U)(\hat\Psi) = 0 \quad\forall \hat\Psi \in \hat X.
\]
The Lagrangian reads
\[
  L(\hat U,\hat Z) = J(\hat U) - \hat A(\hat U)(\hat Z).
\]
As usual, the solution $\hat U$ of the above optimization problem
is characterized by the existence of an adjoint (or dual) solution
$\hat Z \in \hat X$ that satisfies the first-order necessary conditions
\begin{align*}
  0 = L'_U(\hat U,\hat Z)(\hat\Phi)
  &= J'(\hat U)(\hat\Phi) - \hat A'(\hat U)(\hat\Phi,\hat Z)
  &\forall \hat\Phi &\in \hat X, \\
  0 = L'_Z(\hat U,\hat Z)(\hat\Psi)
  &= -\hat A(\hat U)(\hat\Psi)
  &\forall \hat\Psi &\in \hat X.
\end{align*}
The second condition is the given primal problem.
The first condition is the so-called adjoint problem.
The derivatives ${\hat A'}(\hat U)$ and $J'(\hat U)$ are understood in the
Fr\'echet-sense of the nonlinear operator or functional, respectively, evaluated at the
point $\hat U$.

\subsection{Adjoint problem}
From the previous optimality system, we see that the
adjoint problem reads:
find $Z \in X^0$ such that
\begin{equation}\label{Problem: adjoint}
{\hat A}'(\hat U)(\hat\Phi, \hat Z) = J'(\hat U)(\hat\Phi)
\quad \forall \hat\Phi \in \hat X^0.
\end{equation}
We notice that the adjoint problem is always linear, but the nonlinear state
variable $\hU$ enters.

\begin{Problem}[Adjoint problem for stationary fluid-structure interaction]
Find $\hZ \in \hX^0$ such that
\begin{align*}
\hA'(\hU)(\hPsi, \hZ)
  &= \hrho_f \bigl( \hnabla\hpsi^v \hJ\hF^{-1} \hv_f
    + \hnabla\hv_f \hJ\hF^{-1} \hpsi^v , \hz_f^v \bigr)_{\hOmega_f}
    + \hrho_f \bigl( \hnabla\hv_f [\hJ\hF^{-1}]' (\hpsi^u )\hv_f,
    \hz_f^v \bigr)_{\hOmega_f} \\
  &+ \bigl( \hrho_f \nu_f ( \hnabla\hpsi^v \hF^{-1}
    + \hF^{-T} (\hnabla\hpsi^v)^T) \hJ\hF^{-T} ,
    \hnabla\hz_f^v \bigr)_{\hOmega_f} \\
  &+ \bigl( \hrho_f \nu_f (\hnabla\hv_f [\hF^{-1}]'(\hpsi^u )
    + [\hF^{-T}]'(\hpsi^u) \hnabla\hv_f^T ) \hJ\hF^{-T} ,
    \hnabla\hz_f^v \bigr)_{\hOmega_f}
    - \bigl( \hpsi^p \hJ\hF^{-T} , \hnabla\hz_f^v  \bigr)_{\hOmega_f} \\
  &+ \bigl( \hrho_f \nu_f (\hnabla\hv_f \hF^{-1}
    + \hF^{-T} \hnabla \hv_f^T) [\hJ\hF^{-T}]'(\hpsi^u) ,
    \hnabla \hz_f^v \bigr)_{\hOmega_f}
    - \bigl( \hp_f [\hJ\hF^{-T}]'(\hpsi^u) , \hnabla \hz_f^v \bigr)_{\hOmega_f} \\
  &+ \bigl( \lambda_s (\tr\hE'(\hpsi^u)\hF + \tr\hE \hF'(\hpsi^u ))
    + 2\mu_s (\hF'(\hpsi^u )\hE + \hF\hE'(\hpsi^u )) ,
    \hnabla\hz_s^v \bigr)_{\hOmega_s} \\
  &- \bigl( \hpsi^v , \hz_s^u \bigr)_{\hOmega_s}
    + (\alpha_u \hnabla\hpsi^u , \hnabla \hz_f^u )_{\hOmega_f}
    + \bigl( \hat\partial_1 \hpsi^{v_1} + \hat\partial_2 \hpsi^{v_2} ,
    \hz_f^p \bigr)_{\hOmega_f} \\
  &+ \bigl( \hat\partial_2 \hpsi^{u_1} \hat\partial_1 \hv_{f,1}
    - \hat\partial_2 \hpsi^{u_2} \hat\partial_1 \hv_{f,2}
    - \hat\partial_1 \hpsi^{u_2} \hat\partial_2 \hv_{f,1}
    + \hat\partial_1 \hpsi^{u_1} \hat\partial_{f,2} \hv_2 ,
    \hz_f^p \bigr)_{\hOmega_f} \\
    % &\quad + \bigl( \hpsi^p , \hz_s^p \bigr)_{\hOmega_s} \\
  &= J'(\hU)(\hPsi) \qquad \forall \hPsi\in \hX^0.
\end{align*}
\end{Problem}

\subsection{Discrete adjoint problem}
\label{sec_discrete_adjoint}
The corresponding discrete adjoint problem reads
as follows: find $Z_h \in X_h$ corresponding to $U_h \in X_h$ such that
\begin{equation}\label{Problem: discrete adjoint on small space}
{\hat A}'(\hat U_h)(\hat\Phi_h,\hat Z_h)=J'(\hat U_h)(\hat\Phi_h)
\quad \forall \hat\Phi_h \in \hat X_h,
\end{equation}
where $\hat U_h$ solves the discrete primal problem.
Clearly the solution $\hat Z_h$ depends on $\hat U_h$. Algorithmically,
we will first compute the primal solution $\hat U_h$ and then enter
into the adjoint problem \cite{BeRa01}.
Note that due to Galerkin orthogonality, the error estimator stated below
will always yield a zero error unless the adjoint solution
contains higher-order information (see \cite{BeRa01}).
Here, we simply use a globally higher-order finite element solution.
To this end, the adjoint velocities and displacements are approximated
by $Q_c^4$ functions and the pressures by $Q_c^2$ functions.

\subsection{A posteriori error estimators}
In the previous sections, we explained the derivation
and discretization of both the primal and adjoint problems.
Using the main theorem from \cite{BeRa01}, we obtain
the following error identity, which serves as basis
for a posteriori error control.
\begin{theorem}
We have the error identity
\begin{align} \label{dwr_error_representation}
J(\hU) - J(\hU_h) = \frac{1}{2}\rho(\hU_h)(\hZ-\hPhi_h) +
\frac{1}{2}\rho^*(\hU_h, \hZ_h)(\hU-\hPsi_h) + {\cal R}^{(3)}_h
\quad\forall \{\hPsi_h , \hPhi_h\}\in \hX_h \times \hX_h,
\end{align}
where the first two terms on the right-hand side
are given by the primal and adjoint residuals:
\begin{align*}
\rho(\hU_h)(\hZ-\hPhi_h) &:= -A(\hU_h)(\hZ-\hPhi_h), \\ %+ \hat F(\hZ-\hPhi_h) , \\
\rho^*(\hU_h, \hZ_h)(\hU-\hPsi_h) &:= J'(\hU_h)(\hU-\hPsi_h)
- A'(\hU_h)(\hU-\hPsi_h, \hZ_h). % + \hat F(\hU-\hPsi_h).
\end{align*}
The remainder term ${\cal R}^{(3)}_h$ is of cubic order. The arguments
$(\hZ-\hPhi_h)$ and $(\hU-\hPsi_h)$ can be obtained by interpolation
differences, i.e., $(\hZ- i_a \hat Z)$ and $(\hU- i_p \hat U)$, respectively,
where $i_a\colon\hat X_h\high\to \hat X_h$ and $i_p\colon\hat X_h\high\to \hat X_h$
are interpolations from higher order finite element spaces to low-order
spaces; see also Section \ref{sec_discrete_adjoint}.
\end{theorem}
\begin{proof}
We refer the reader to \cite{BeRa01,RanVi2013,EndtLaWi20}.
\end{proof}
\begin{Corollary}[Error estimator based on primal residual]
\label{coro_1}
The primal error identity reads:
\begin{align} \label{dwr_error_representation_primal}
J(\hU) - J(\hU_h) = \rho(\hU_h)(\hZ-\hPhi_h) + {\cal R}^{(2)}_h,
\end{align}
where the remainder term is now of second order.
\end{Corollary}
The previous error identity can be used to define the error estimator $\eta$,
which can be further utilized to design adaptive schemes.
Therein, we propose practical error estimators in which all
information can be computed and remainder terms are neglected.

\begin{Definition}[Practical error estimators]
A practical error estimator for the goal functional $J(\hat U)$ reads:
\begin{equation}
\label{eq_err_est_both}
\eta_h:=\frac{1}{2}\rho(\hat U_h)(\hat Z_h - i_a \hat Z_h)
+ \frac{1}{2}\rho^*(\hat U_h, \hat Z_h)(\hat U_h - i_p \hat U_h).
\end{equation}
A purely primal and hence less accurate, but cheaper,
practical error estimator reads:
\begin{equation}
\label{eq_err_est_prim}
\eta_h:= \frac{1}{2}\rho(\hat U_h)(\hat Z_h - i_a \hat Z_h).
\end{equation}
The primal error part has a remainder term of only second order
(see Corollary \ref{coro_1}), but is cheaper because only the adjoint
solution needs to contain higher order information for the interpolation
$i_a$. In \eqref{eq_err_est_both} both
the primal and the adjoint must contain higher order information for
$i_a$ and $i_p$, which is
of course expensive in terms of the computational cost.
\end{Definition}

\subsection{Multiple goal functionals}

The error estimation is now extended to multiple goal functionals.
Assume that we are given $N$ goal functionals $J_i$, $i=1,\dots,N$.
In a flow problem this might be drag and lift as well as
estimates of the pressure and some solid displacement
of the elastic structure. In this case, we would have $N=4$.

We construct an overall goal functional $J_c$ as a
convex combination of the individual goal functionals,
with weights $\omega_i \ge 0$ that sum up to one
and signs $\sigma_i \in \{-1, +1\}$:
\begin{equation}\label{definition:J_c_tilde}
  J_c(\hat\Phi) :=
  \sum_{i=1}^{N} \omega_i \sigma_i J_i(\hat\Phi),
  \quad \hat\Phi \in \hat X.
\end{equation}
\iffalse
If the weights are similar the relative errors of all
single functionals are of the same order.
\fi
Here the choice of $\sigma_i$ is a crucial aspect
since all terms in the sum should have the same sign to avoid cancellation.
As we need to compute $|J_c(\hat U) - J_c(\hat U_h)|$,
we follow previous studies and use
\begin{equation}\label{weight w_i}
\sigma_i := \text{sign}(J_i(\hat U) - J_i(\hat U_h)).
\end{equation}
Later, in Section \ref{sec_tests}, we set
\[
w_i := \omega_i \sigma_i.
\]
We notice that a relative combined functional may be defined with 
$w_i := \omega_i \frac{\sigma_i}{|J_i(\hat U_h)|}$, which is however not employed in 
this work. Some computations studying the influence of different 
weights $\omega_i$ are conducted in Section \ref{sec_tests}.

A first sign computation was proposed in \cite{HaHou03} and later extended
to a more efficient way in \cite{EndtWi17}.
In practice $J_c$ from \eqref{definition:J_c_tilde} is now
used as right hand side in the adjoint problems \eqref{Problem: adjoint}
and \eqref{Problem: discrete adjoint on small space}, respectively.
We notice that for nonlinear goal functionals, we need
the Fr\'echet derivative of the combined functional $J_c$ \eqref{definition:J_c_tilde}; for
technical details, we refer to \cite{EndtLaWi18}.

\begin{Definition}[Combined functional error estimators]
A practical error estimator for the combined goal functional $J_c(\hat U)$ reads:
\begin{equation}
\label{eq_err_est_both_comb}
J_c(\hat U) - J_c(\hat U_h) \approx \eta_h:=\frac{1}{2}\rho(\hat U_h)(\hat Z_h - i_a \hat Z_h)
+ \frac{1}{2}\rho^*(\hat U_h, \hat Z_h)(\hat U_h - i_p \hat U_h).
\end{equation}
As above, a purely primal and less accurate but cheaper
practical error estimator reads:
\begin{equation}
\label{eq_err_est_prim_comb}
J_c(\hat U) - J_c(\hat U_h) \approx \eta_h:= \frac{1}{2}\rho(\hat U_h)(\hat Z_h - i_a \hat Z_h).
\end{equation}
In both cases the primal problem is computed as for a single
goal functional, and in the adjoint problem the combined
goal functional $J_c$ (including sign computation) is employed.
\end{Definition}

\section{Error localization and adaptive algorithms}
\label{sec_error_loc}

\subsection{Localization}
For the localization of the error estimator to single DoFs or elements,
we use the partition-of-unity (PU) technique as suggested in \cite{RiWi15_dwr}.
To this end, we choose a set of finite element basis functions
$\hV_{PU}:= \{\hat \psi_1, \dots,\hat\psi_M\}$
with $\dim\hat V_{PU}=M$ such that
$\sum_{i=1}^{N} \psi_i \equiv 1$. These functions can be low-order
scalar-valued bilinear shape functions $Q_1^c$.
We then distribute $\eta_i$ to the corresponding elements with certain weights as explained in \cite{EndtLaWi18,Endt21}.
Inserting this into \eqref{eq_err_est_prim_comb},
we obtain
\begin{Proposition}
\label{prop_41}
For the combined goal functional $J_c$, using the primal error part
$\rho(\hU_h)(\cdot)$, we have the a posteriori error estimate
 \begin{equation}
\label{eq_est}
|J_c(\hU) - J_c(\hU_h)| \leq |\eta_h| := \bigl| \sum_{i=1}^M \eta_i \bigr| \leq \sum_{i=1}^M |\eta_i|
 \end{equation}
with the PU-DoF indicators
\begin{align*}
\eta_i &= - A(\hU_h)((\hZ_h - i_a \hZ_h)\hPsi_i)
%+ \hat F((\hZ_h^{(2)}-i_h \hZ_h^{(2)})\hPsi_i)
\\
&= -
(\hrho_f \hJ  (\hF^{-1}\hv_f\cdot\hnabla) \hv_f),\hpsi_i^v)_{\hOmega_f}
- (\hJ\hsigma_f\hF^{-T},\hnabla\hpsi_i^v)_{\hOmega_f}
+ \langle \hg_f, \hpsi_i^v \rangle_{\hGamma_N}\\
%+ (\hrho_f \hJ\hf_f, \hpsi^v)_{\hOmega_f}\\
&\quad - (\hF\hSigma ,\hnabla\hpsi_i^v)_{\hOmega_s}
%+ (\hrho_s\hf_s, \hpsi^v)_{\hOmega_s}
- (\hsigma\tsb{mesh} , \hnabla\hpsi_i^u)_{\hOmega_f}
- (\hdiv\,(\hJ\hF^{-1}\hv_f),\hpsi_i^p)_{\hOmega_f}
+ (\hJ\hf_f,\hpsi_i^v) + (\hf_s,\hpsi_i^v).
\end{align*}
Here the \textit{weighting functions} are defined
with the interpolation $i_a\colon\hX_h\high\to\hX_h$ as
\begin{align*}
\hpsi_i^v &:= (\hat z_{h,v}\high-\hat z_{h,v})\psi_i, &
\hpsi_i^u &:= (\hat z_{h,u}\high-\hat z_{h,u})\psi_i, &
\hpsi_i^p &:= (\hat z_{h,p}\high-\hat z_{h,p})\psi_i.
\end{align*}
\end{Proposition}
To measure the quality of the proposed error estimator,
we consult the effectivity index \cite{BaRhei78b} and indicator index \cite{RiWi15_dwr},
respectively:
\begin{align}
  \Ieff &:= \frac{|\eta_h|}{|J_c(\hat U)-J_c(\hat U_h)|}, &
  \Iind &:= \frac{\sum_i|\eta_i|}{|J_c(\hat U)-J_c(\hat U_h)|}.
\end{align}

\subsection{Adaptive algorithm}
\begin{enumerate}
\item Compute the primal solution $\hU_h\in \hat X_h$ on the mesh ${\cal T}_l$,
  where $l\in\mathbb{N}$ is the current mesh level.
%\item To construct $J_c$, first determine the signs $\sigma_i$
%  via \eqref{weight w_i}.
\item Construct the combined goal functional $J_c$
  via \eqref{definition:J_c_tilde}.
\item Solve \eqref{Problem: discrete adjoint on small space}
  on the mesh ${\cal T}_l$ with $J_c$ as right hand side
  and obtain (high-order) $\hZ_h\in \hat X_h\high$.
\item Evaluate $|\eta| := |\sum_{i} \eta_{i}|$ in \eqref{eq_est}.
\item If the stopping criterion is satisfied,
$|J_c(\hU) - J_c(\hU_h)| \leq |\eta| \leq TOL$, then accept $U_h$
within the tolerance $TOL$.
Otherwise, proceed
to the following step.
\item Mark all elements $K_i$ for refinement that touch DoFs $i$
whose indicator $\eta_{i}$ satisfies
$\eta_{i} \geq \frac{\alpha\eta}{M\tsb{el}}$ (where
$M\tsb{el}$ denotes the total number of elements of the mesh $\mathcal{T}_h$
and $\alpha \approx 1$).
\item Refine all marked elements to obtain the mesh ${\cal T}_{l+1}$.
\item Go to Step 1.
\end{enumerate}

\section{Numerical tests}
\label{sec_tests}
In this section, we present three numerical tests:
an elastic lid-driven cavity with two goal functionals,
an elastic bar in a chamber with two goal functionals, and the FSI-1 benchmark
with three goal functionals.
The implementation
is based on the open-source finite element library deal.II \cite{dealII91,deal2020}
and extensions of our own fluid-structure code publications \cite{Wi13_fsi_with_deal,Wi21_git}
towards multigoal-oriented error estimation from \cite{Endt21}. In all examples,
we work with the primal error estimator from Proposition \ref{prop_41}.

%%%%%%%%%%%%%%%%%%%%%%%%%%%%%%%%%%%%%%%%%%%%%%%%%%%%%%%%%%%
\subsection{Example 1: lid-driven cavity with elastic volume at bottom}
This first configuration is a well-known
example in computational fluid dynamics and was extended
to fluid-structure interaction in \cite{Du06} where an elastic bottom
is added in the lower part of the domain.%

\paragraph{Configuration}
The computational domain is $\Omega = (0,2)^2$ with flow
in $\Omega_f = [0, 2] \times [0.5, 2]$
and the solid domain $\Omega_s = [0, 2] \times [0, 0.5]$.

\paragraph{Boundary conditions}
On the top boundary $\hGamma\tsb{top} = \{2\} \times [0, 2]$,
we prescribe overflow:
\[
v_0 = 0.5 \times
\begin{cases}
\sin^2 (\pi x/0.6), & x\in [0.0,0.3],\\
1, & x\in (0.3,1.7),\\
\sin^2 (\pi (x-2.0)/0.6), & x\in [1.7,2.0].
\end{cases}
\]
On the other boundaries, we use homogeneous Dirichlet conditions, i.e., $\hv=0$ on
$\partial\hOmega\setminus\hGamma\tsb{top}$. Moreover, $\hu=0$ on $\partial\hOmega$
and $\partial_n p = 0$ on $\partial\hOmega$. We notice that, because of the pressure
boundary conditions, the pressure is not unique and must be constrained.

\paragraph{Parameters}

For the fluid
we use the density $\varrho_f = \SI{1.0}{kg m^{-3}}$ and kinematic
viscosity $\nu_f = \SI{0.2}{m^2 s^{-1}}$. The elastic
solid is characterized by
the density $\varrho_s = \SI{1.0}{kg m^{-3}}$ and Poisson's ratio $\nu_s = 0.4$. Furthermore,
we use the Lam\'e coefficient $\mu_s = \SI{2.0}{kg m^{-1} s^{-2}}$.

\paragraph{Quantities of interest and goals}
In this example, we construct the following combined functional:
\[
J_c(\hU) = w_1 \Jdrag(\hU) + w_2 J_2(\hu),
\]
with the weights $w_1,w_2\in\mathbb{R}$ and $\omega_1=\omega_2 = 0.5$
(see \eqref{definition:J_c_tilde}) and $\Jdrag(\hU)$ (see definition
below \eqref{drag_lift_forces}) and
$J_2(\hu) := \hu(1.5,0.25)$.
The individual reference values are computed on a sufficiently refined mesh:
\iffalse
\[
  J_c(\hat U) = \num[round-mode=places,round-precision=5]{-4.91167835466135e-02}.
\]
\fi
\begin{align*}
  \Jdrag(\hU)
  &= \num[round-mode=places,round-precision=5]{-9.3543731705807223e-02}, \\
  J_2(\hU)
  &= \num[round-mode=places,round-precision=5]{-4.6898353874198270e-03}.
\end{align*}
The combined reference value is then obtained as
\begin{equation}
\label{eq_ex_1_Jc}
J_c(\hat U) = \omega_1 \Jdrag(\hU) + \omega_2 J_2(\hu).
\end{equation}

\paragraph{Results and discussion}
The original geometry is once uniformly refined and serves as initial mesh.
Then we perform six adaptive refinements.
The final two goal functionals and the combined functional
values are provided in Table \ref{tab:ex1-goal-final}.
We see that the values do not necessarily need to be positive.

\begin{table}[ht]
  \centering
  \sisetup{round-mode=places,round-precision=3,table-format=+1.3e+1}
  \begin{tabular}{lS}
    \toprule
    Drag     &  -9.35496118e-02 \\
    Pressure &  -4.68801063e-03 \\
    $J_c$    &  -4.91188112e-02 \\
     \bottomrule
  \end{tabular}
  \caption{Example 1: Final values of goal functionals for $\omega_1 = \omega_2 = 0.5$,
where $J_c$ is obtained with the help of \eqref{eq_ex_1_Jc}.}
  \label{tab:ex1-goal-final}
\end{table}

As we observe in
Table \ref{tab:ex1-levels},
the estimated error $\eta_h$ of $J_c$ drops down to \num{1.71e-06}.
The effectivity indices $\Ieff$ perform relatively well being close to $1$
on the last three meshes. The slight deviation can be explained by the
nonlinear, coupled problem and by utilizing only the primal error part.
In a comparison with uniform mesh refinement, we observe that the true
error behaves first sightly worse than adaptive mesh refinement, but finally better.
However the computational cost measured in terms of the degrees of freedom
is much higher for uniform mesh refinement. The estimators behave similarly
at each refinement level, but for different degrees of freedom. These results
clearly show that adaptive mesh refinement is an efficient procedure.
Table~\ref{tab:ex1-levels_b} shows the effect of choosing different
weights $\omega_1, \omega_2$. Therein, always the above reference
value for $J_c(\hat U)$ is used.
For the case $\omega_1= 0.00$ and $\omega_2=1.00$ all
values vanish. The reason is that our point evaluation is a grid
point and therefore the values vanish if we use interpolation.
A more detailed explanation for this phenomena and solutions to
overcome this problem are found in \cite{EnLaWi21a}.

\begin{table}[ht]
  \centering
  \sisetup{table-number-alignment=center,table-format=1.2e+1}
  \begin{tabular}{rSSS
    %S[scientific-notation=fixed,fixed-exponent=0,table-format=1.3]
    %S[scientific-notation=fixed,fixed-exponent=0,table-format=2.2]}
    SS}
      \toprule
    Dofs&${|J_c(\hat U)-J_c(\hat U_h)|}$&{$|\eta_h|$}&{$\sum_i|\eta_i|$}&$\Ieff$&$\Iind$ \\
    \midrule
%    195   & 1.76e-03 & 3.62e-03 & 1.99e-02 & 2.06e+00 & 1.13e+01 \\
%    554   & 1.85e-03 & 2.22e-03 & 5.47e-03 & 1.20e+00 & 2.96e+00 \\
%    1671  & 6.93e-04 & 3.50e-04 & 2.42e-03 & 5.05e-01 & 3.50e+00 \\
%    3848  & 6.27e-05 & 7.39e-05 & 9.56e-04 & 1.18e+00 & 1.53e+01 \\
%    9925  & 5.61e-06 & 7.98e-06 & 2.66e-04 & 1.42e+00 & 4.75e+01 \\
%    20687 & 4.06e-06 & 3.41e-06 & 7.42e-05 & 8.42e-01 & 1.83e+01 \\
%195   &  1.76e-03 &       1.81e-03  &      9.96e-03 &       1.03e+00 &       5.66e+00\\
%554   &  1.85e-03 &       1.11e-03  &      2.74e-03 &       6.02e-01 &       1.48e+00\\
%1671  &  6.93e-04 &       1.75e-04  &      1.21e-03 &       2.53e-01 &       1.75e+00\\
%3848  &  6.27e-05 &       3.69e-05  &      4.78e-04 &       5.89e-01 &       7.63e+00\\
%9925  &  5.61e-06 &       3.99e-06  &      1.33e-04 &       7.11e-01 &       2.38e+01\\
%20687 &  4.06e-06 &       1.71e-06  &      3.71e-05 &       4.21e-01 &       9.14e+00\\ \midrule
195	&	8.80e-04&	1.81e-03&	9.96e-03&	2.06e+00&	1.13e+01\\
554	&	9.23e-04&	1.11e-03&	2.74e-03&	1.20e+00&	2.96e+00\\
1671	&	3.46e-04&	1.75e-04&	1.21e-03&	5.05e-01&	3.50e+00\\
3848	&	3.14e-05&	3.69e-05&	4.78e-04&	1.18e+00&	1.53e+01\\
9925	&	2.80e-06&	3.99e-06&	1.33e-04&	1.42e+00&	4.75e+01\\
20687	&	2.03e-06&	1.71e-06&	3.71e-05&	8.42e-01&	1.83e+01\\ \midrule
    \multicolumn6c{Uniform mesh refinement} \\
    \midrule
%195   &  1.76e-03  &      1.81e-03  &      9.96e-03 &       1.03e+00  &      5.66e+00\\
%657   &  1.63e-03  &      1.06e-03  &      2.75e-03 &       6.48e-01  &      1.68e+00\\
%2397  &  4.52e-04  &      1.55e-04  &      1.16e-03 &       3.44e-01  &      2.58e+00\\
%9141  &  5.57e-05  &      1.52e-05  &      4.05e-04 &       2.73e-01  &      7.27e+00\\
%35685 &  6.31e-07  &      2.55e-06  &      1.17e-04 &       4.05e+00  &      1.85e+02\\
%140997&  9.80e-07  &      1.52e-06  &      3.22e-05 &       1.55e+00  &      3.29e+01\\ \midrule
195	&	8.80e-04&	1.81e-03&	9.96e-03&	2.06e+00&	1.13e+01\\
657	&	8.15e-04&	1.06e-03&	2.75e-03&	1.30e+00&	3.37e+00\\
2397	&	2.26e-04&	1.55e-04&	1.16e-03&	6.88e-01&	5.15e+00\\
9141	&	2.78e-05&	1.52e-05&	4.05e-04&	5.46e-01&	1.45e+01\\
35685	&	3.16e-07&	2.55e-06&	1.17e-04&	8.09e+00&	3.70e+02\\
140997	&	4.90e-07&	1.52e-06&	3.22e-05&	3.10e+00&	6.57e+01\\
    \bottomrule
  \end{tabular}
  \caption{Example 1: Degrees of freedom, true error, estimator and indices.}
  \label{tab:ex1-levels}
\end{table}

\begin{table}[ht]
  \centering
  \sisetup{table-number-alignment=center,table-format=1.2e+1}
  \begin{tabular}{rrrSSS
    %S[scientific-notation=fixed,fixed-exponent=0,table-format=1.3]
    %S[scientific-notation=fixed,fixed-exponent=0,table-format=2.2]}
    SS}
      \toprule
$\omega_1$&$\omega_2$&Dofs&${|J_c(\hat U)-J_c(\hat U_h)|}$&{$|\eta_h|$}&{$\sum_i|\eta_i|$}&$\Ieff$&$\Iind$ \\
    \midrule
0.00 &1.00  &  195      &       6.20e-04&	0.00e+00&	0.00e+00&	0.00e+00&	0.00e+00\\ \midrule
0.25 &0.75  &  195	&       1.30e-04&	9.04e-04&	4.98e-03&	6.96e+00&	3.84e+01\\
0.25 &0.75  &  554	&       6.18e-04&	5.56e-04&	1.37e-03&	8.99e-01&	2.21e+00\\
0.25 &0.75  &  1671	&       3.19e-04&	8.75e-05&	6.05e-04&	2.74e-01&	1.89e+00\\
0.25 &0.75  &  3848	&       1.89e-05&	2.45e-05&	2.96e-04&	1.29e+00&	1.56e+01\\
0.25 &0.75  &  9448	&       5.03e-06&	5.87e-06&	7.91e-05&	1.17e+00&	1.57e+01\\
0.25 &0.75  &  23303	&       3.00e-07&	1.84e-06&	2.21e-05&	6.12e+00&	7.36e+01\\ \midrule
0.75 &0.25  &  195	&       1.63e-03&	2.71e-03&	1.49e-02&	1.66e+00&	9.17e+00\\
0.75 &0.25  &  554	&       1.23e-03&	1.67e-03&	4.10e-03&	1.36e+00&	3.34e+00\\
0.75 &0.25  &  1671	&       3.73e-04&	2.63e-04&	1.82e-03&	7.04e-01&	4.87e+00\\
0.75 &0.25  &  3848	&       4.38e-05&	4.94e-05&	6.62e-04&	1.13e+00&	1.51e+01\\
0.75 &0.25  &  9865	&       1.07e-05&	1.38e-06&	1.89e-04&	1.29e-01&	1.77e+01\\
0.75 &0.25  &  20512	&       1.15e-06&	1.07e-06&	5.33e-05&	9.32e-01&	4.63e+01\\ \midrule
1.00 &0.00  &  195	&       2.38e-03&	3.62e-03&	1.99e-02&	1.52e+00&	8.37e+00\\
1.00 &0.00  &  554	&       1.53e-03&	2.22e-03&	5.47e-03&	1.45e+00&	3.57e+00\\
1.00 &0.00  &  1671	&       4.00e-04&	3.50e-04&	2.42e-03&	8.75e-01&	6.05e+00\\
1.00 &0.00  &  3848	&       5.62e-05&	6.18e-05&	8.61e-04&	1.10e+00&	1.53e+01\\
1.00 &0.00  &  9797	&       1.66e-06&	1.31e-06&	2.46e-04&	7.88e-01&	1.48e+02\\
1.00 &0.00  &  18985	&       2.29e-06&	5.09e-07&	7.11e-05&	2.23e-01&	3.11e+01\\
    \bottomrule
  \end{tabular}
  \caption{Example 1: Degrees of freedom, true error, estimator and indices
    for different weights $\omega_1, \omega_2$.
For each weight combination, the reference $J_c$ is obtained 
with the help of \eqref{eq_ex_1_Jc} and the values for 
$\Jdrag(\hU)$ and $J_2(\hu)$ listed in Table \ref{tab:ex1-goal-final}.
%
%    \MCS{and the reference value $J_c(\hat U) =
%      \num[round-mode=places,round-precision=5]{-9.8233567093227045e-02}$.
%      We need different reference values. How do we present them?
%      More importantly:
%      are the results here based on the single given reference value?
    }
  \label{tab:ex1-levels_b}
\end{table}

\newpage
In Table \ref{tab:ex1-levels_c}, the values of the individual goal functionals
and the combined goal functional are listed.
These allow us to study the evolution
of the individual goal functionals under mesh refinement for the different weights.
For each weight combination, the reference $J_c$ is obtained 
with the help of \eqref{eq_ex_1_Jc} and the values for 
$\Jdrag(\hU)$ and $J_2(\hu)$ listed in Table \ref{tab:ex1-goal-final}.
This study is important insofar as we only control $J_c$ via the adjoint problem
in our multigoal framework, but not directly the individual goals.

\begin{table}[ht]
  \centering
  \sisetup{round-mode=places,round-precision=2,
    table-number-alignment=center,table-format=+1.2e+1}
  \begin{tabular}{rrrSSS}
      \toprule
    $\omega_1$&$\omega_2$&Dofs&{$J_c(\hat U)$}&{$\Jdrag(\hU)$}&{$J_2(\hU)$} \\
    \midrule
0.00 &1.00  &  195      &    -4.06961914e-03&	-9.59236920e-02&	-4.06961914e-03\\ \midrule
0.25 &0.75  &  195	&    -2.70331373e-02&	-9.59236920e-02&	-4.06961914e-03\\
0.25 &0.75  &  554	&    -2.62851131e-02&	-9.20115888e-02&	-4.37628787e-03\\
0.25 &0.75  &  1671	&    -2.65838614e-02&	-9.31437330e-02&	-4.39723760e-03\\
0.25 &0.75  &  3848	&    -2.68843666e-02&	-9.34875668e-02&	-4.68329984e-03\\
0.25 &0.75  &  9448	&    -2.68982778e-02&	-9.35469997e-02&	-4.68203721e-03\\
0.25 &0.75  &  23303	&    -2.69030092e-02&	-9.35476790e-02&	-4.68811925e-03\\ \midrule
0.75 &0.25  &  195	&    -7.29601738e-02&	-9.59236920e-02&	-4.06961914e-03\\
0.75 &0.25  &  554	&    -7.01027636e-02&	-9.20115888e-02&	-4.37628787e-03\\
0.75 &0.25  &  1671	&    -7.09571091e-02&	-9.31437330e-02&	-4.39723760e-03\\
0.75 &0.25  &  3848	&    -7.12865000e-02&	-9.34875668e-02&	-4.68329984e-03\\
0.75 &0.25  &  9865	&    -7.13195520e-02&	-9.35336861e-02&	-4.67714961e-03\\
0.75 &0.25  &  20512	&    -7.13314072e-02&	-9.35461273e-02&	-4.68724710e-03\\ \midrule
1.00 &0.00  &  195	&    -9.59236920e-02&	-9.59236920e-02&	-4.06961914e-03\\
1.00 &0.00  &  554	&    -9.20115888e-02&	-9.20115888e-02&	-4.37628787e-03\\
1.00 &0.00  &  1671	&    -9.31437330e-02&	-9.31437330e-02&	-4.39723760e-03\\
1.00 &0.00  &  3848	&    -9.34875668e-02&	-9.34875668e-02&	-4.68329984e-03\\
1.00 &0.00  &  9797	&    -9.35420689e-02&	-9.35420689e-02&	-4.68157420e-03\\
1.00 &0.00  &  18985	&    -9.35460179e-02&	-9.35460179e-02&	-4.68650034e-03\\
    \bottomrule
  \end{tabular}
  \caption{Example 1: Actual values for the individual goals and combined goal functional.
}
  \label{tab:ex1-levels_c}
\end{table}

\begin{figure}[ht]
\centering
\plot{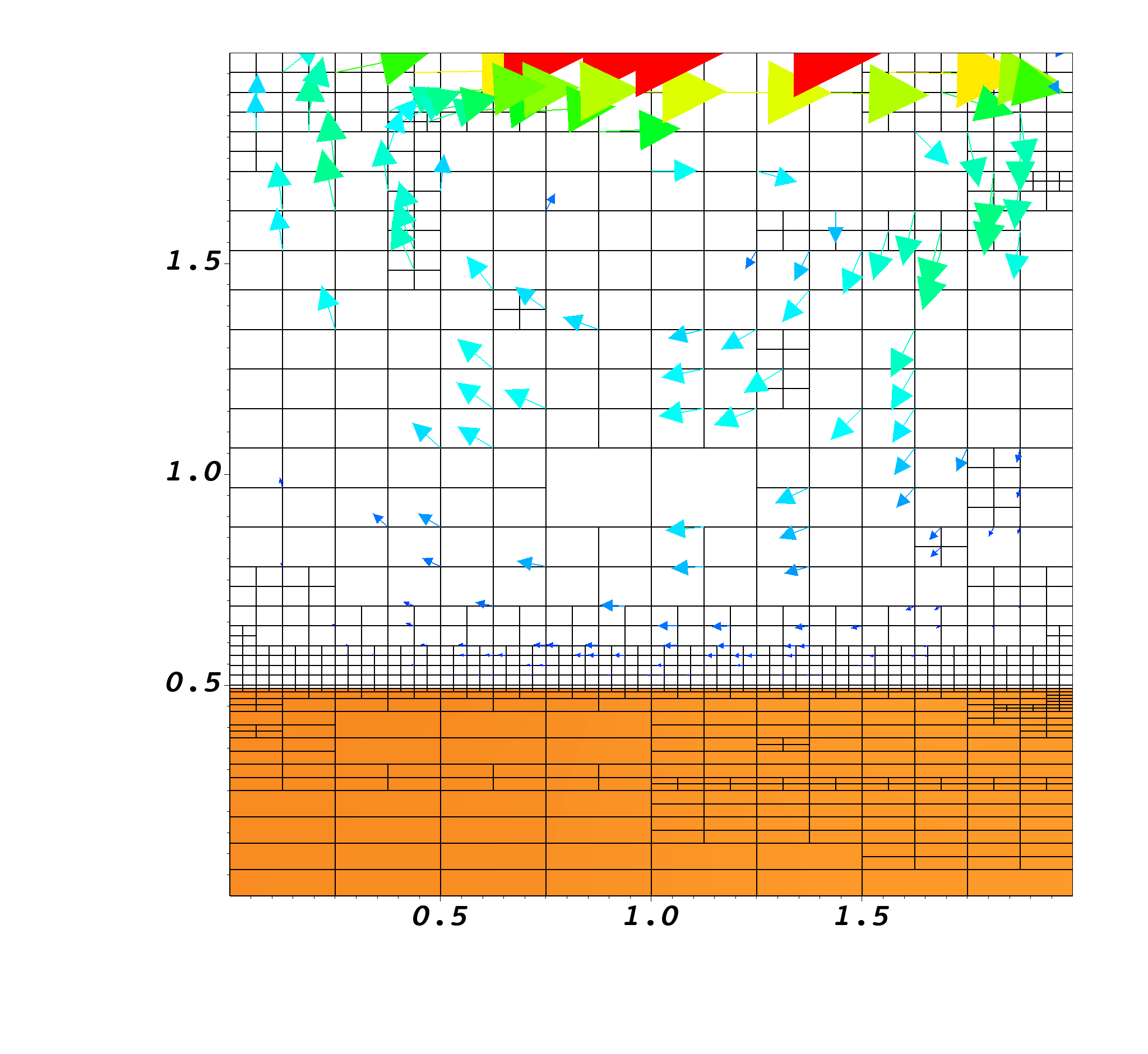}
\plot{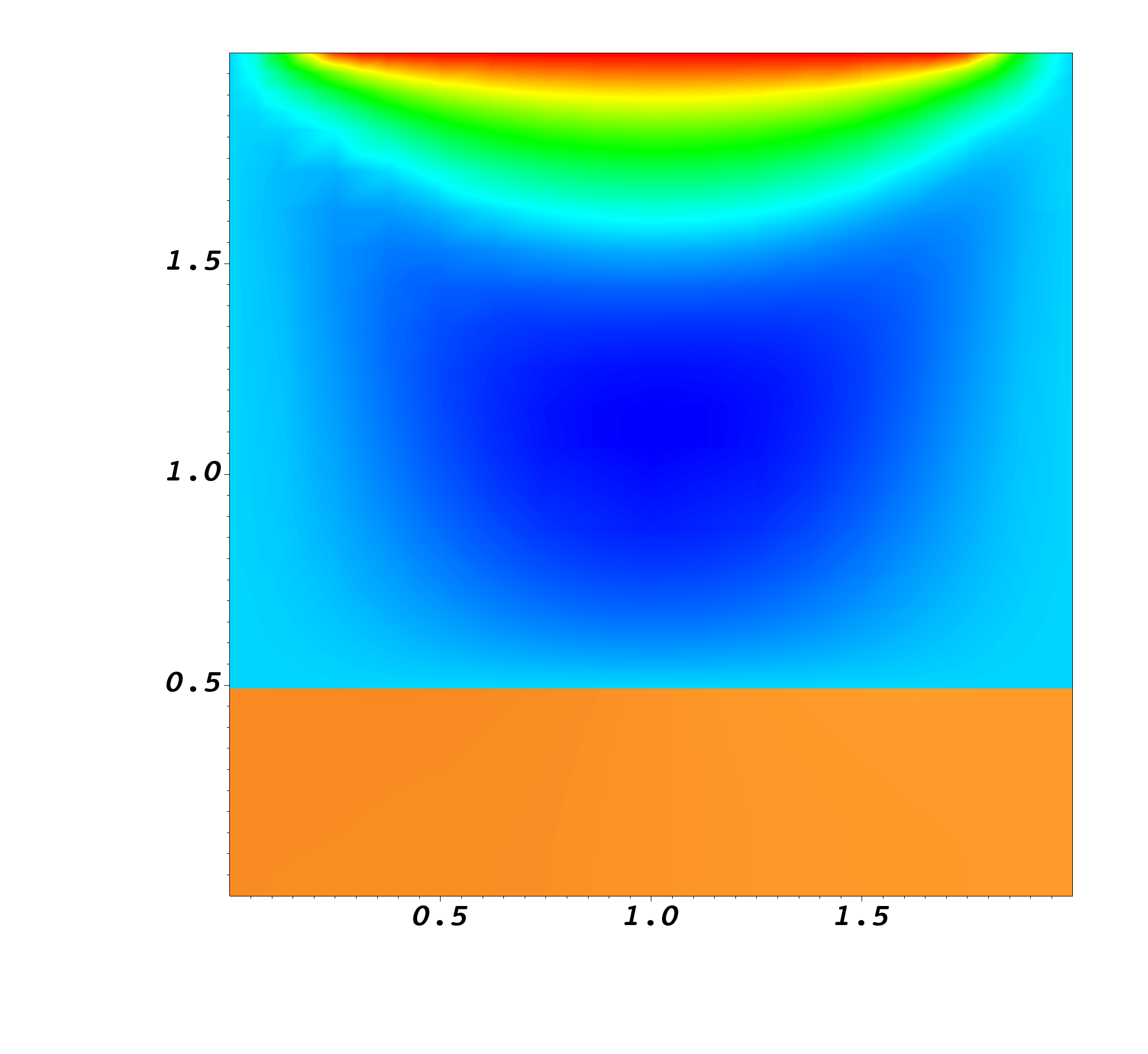}
\caption{Example 1: Left: adaptive mesh and vector plot. We observe adaptive
mesh refinement along the FSI interface and also around the point value
evaluation of $J_2(\hu) := \hu(1.5,0.25)$. Right: pseudocolor plot of $\hv_x$. The lower
brown part shows the elastic solid domain $\hOmega_s$.}
\label{pic_ex_2a}
\end{figure}

\begin{figure}[ht]
\centering
\plot{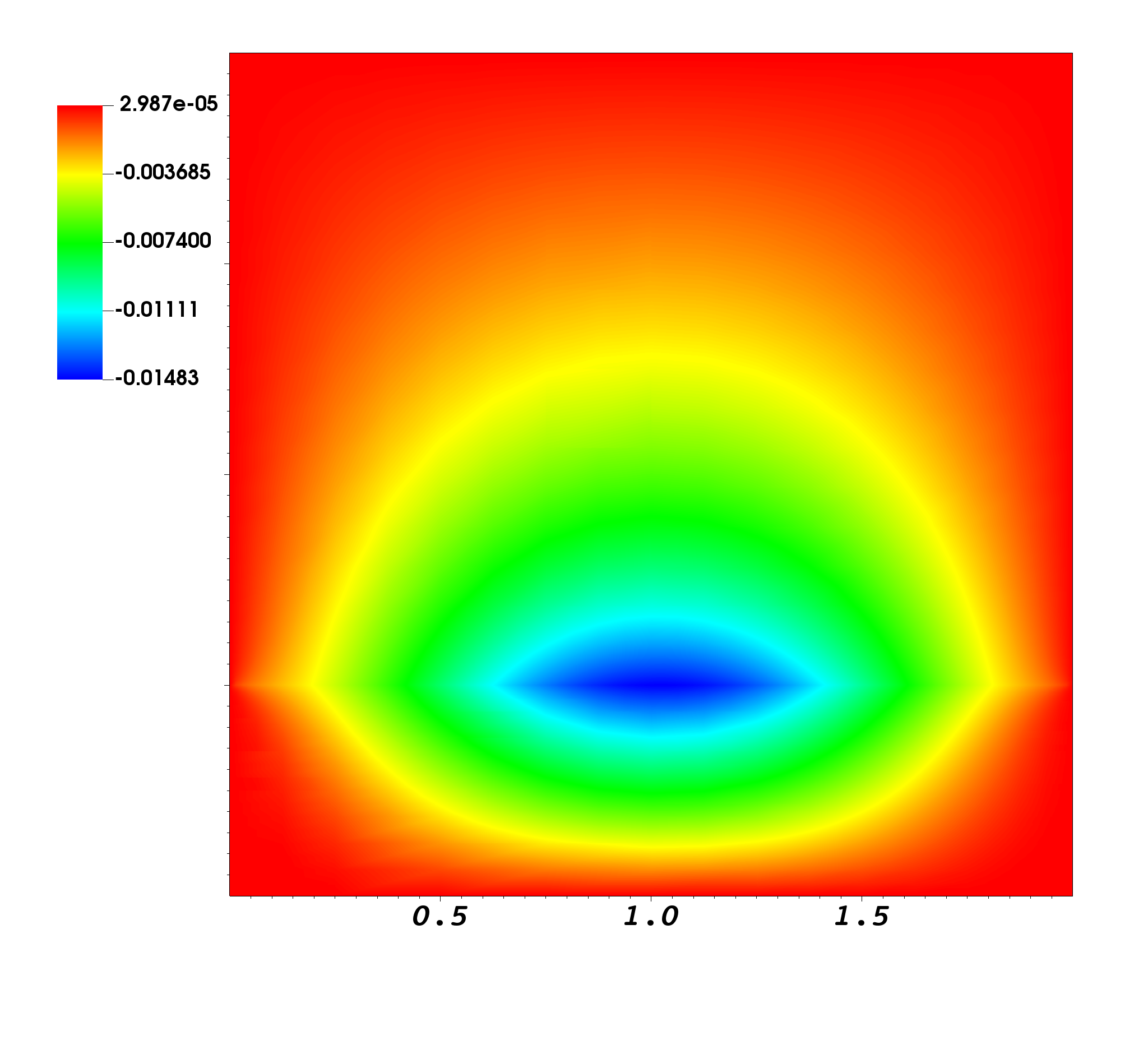}
\plot{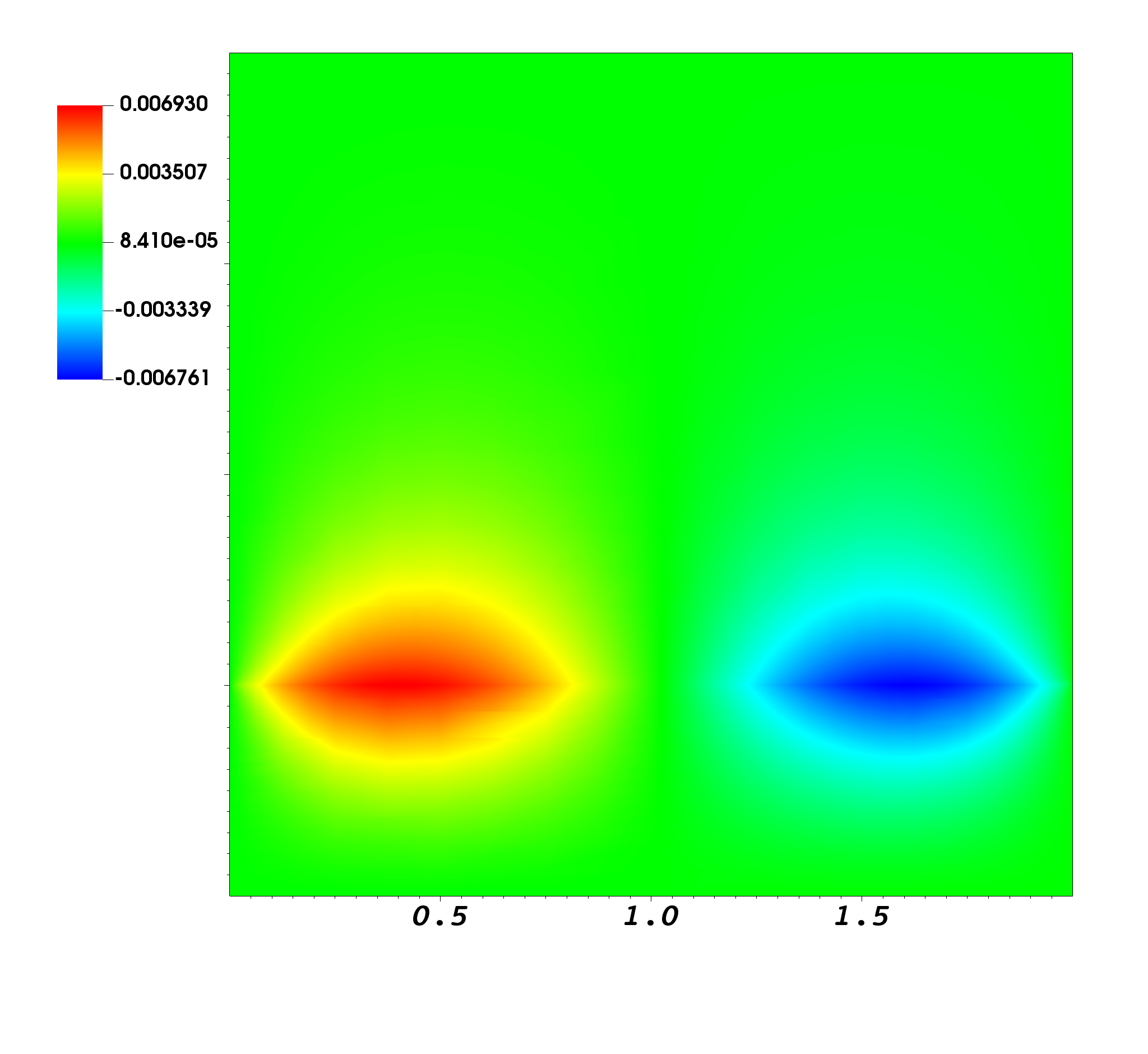}
\caption{Example 1: Left and right: pseudocolor plots of the primal displacement solutions $\hu_x$ and
$\hu_y$, respectively.}
\label{pic_ex_2b}
\end{figure}

\begin{figure}[ht]
\centering
\plot{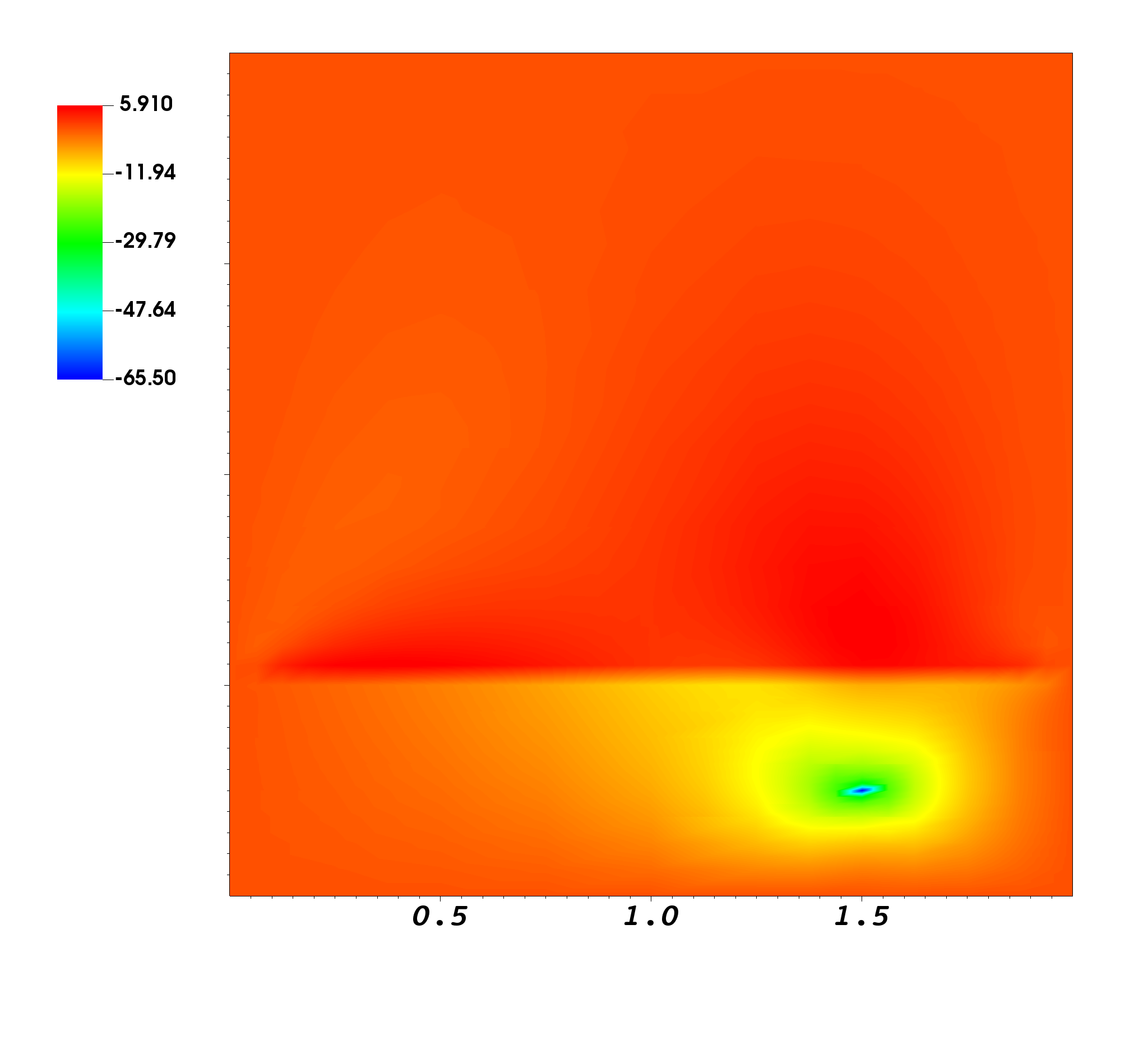}
\plot{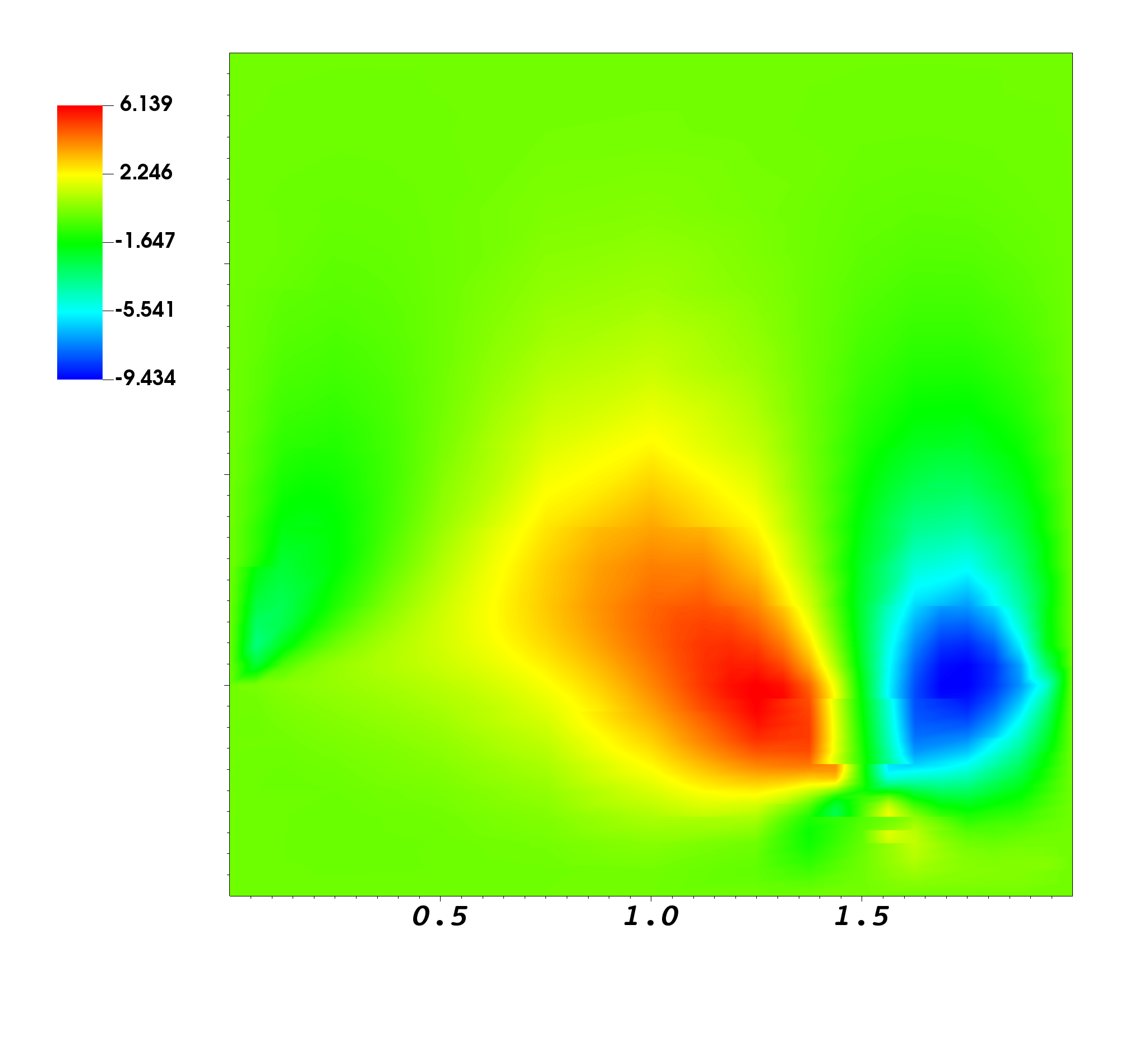}
\caption{Example 1: Left and right: pseudocolor plots of the adjoint solutions $z^{\hv_x}$ and
$z^{\hv_y}$, respectively.}
\label{pic_ex_2c}
\end{figure}

\newpage
In Figure \ref{pic_ex_2a} the primal solution is shown with the adaptively
refined mesh. Refinements with respect to both goal functionals
can be observed. In Figure \ref{pic_ex_2b}, the displacement fields
are displayed.
Afterward, in Figure \ref{pic_ex_2c}, the adjoint solutions
of the velocities are shown. The overflow
velocity boundary condition introduces pressure singularities in both upper corners.
We also notice that at the boundary points where the FSI interface intersects, different
boundary and interface conditions interact, which may also lead to
a slight degeneration of the effectivity indices.

% Internal comments
%Thomas Dunne 2006 hatte $\|v\|^2$ , was okay ist.
%
%Vielleicht koennte man hier folgendes machen:
%
%1. $\|v\|^2$ um mit Thomas Dunne zu vergleichen
%
%2. Dann $\|v\|^2$ , den stress auf $\Gamma_i$ in y-direction, und
%    dann vielleicht ein Point displacement z.B. in $(0.5, 0.5)$.
%    Zur Erinnerung: das Gebiet ist $(0,2)^2$. Und das Interface
%    ist $\Gamma_i = {(x,y) | y = 0.5}$
%

\newpage
%%%%%%%%%%%%%%%%%%%%%%%%%%%%%%%%%%%%%%%%%%%%%%%%%%%%%%%%%%%
\subsection{Example 2: elastic bar in a chamber}
This second configuration is taken from \cite{RiWi13}.
Here, the flow is driven by a pressure difference.
This configuration is challenging for the ALE transformation because
a higher inflow pressure yields higher fluid flow and will
cause the elastic beam to close the small upper channel. Consequently,
both the adjoint equation and the mesh motion technique play important
roles.

\begin{figure}[H]
  \centering
  \includegraphics[width=78mm]{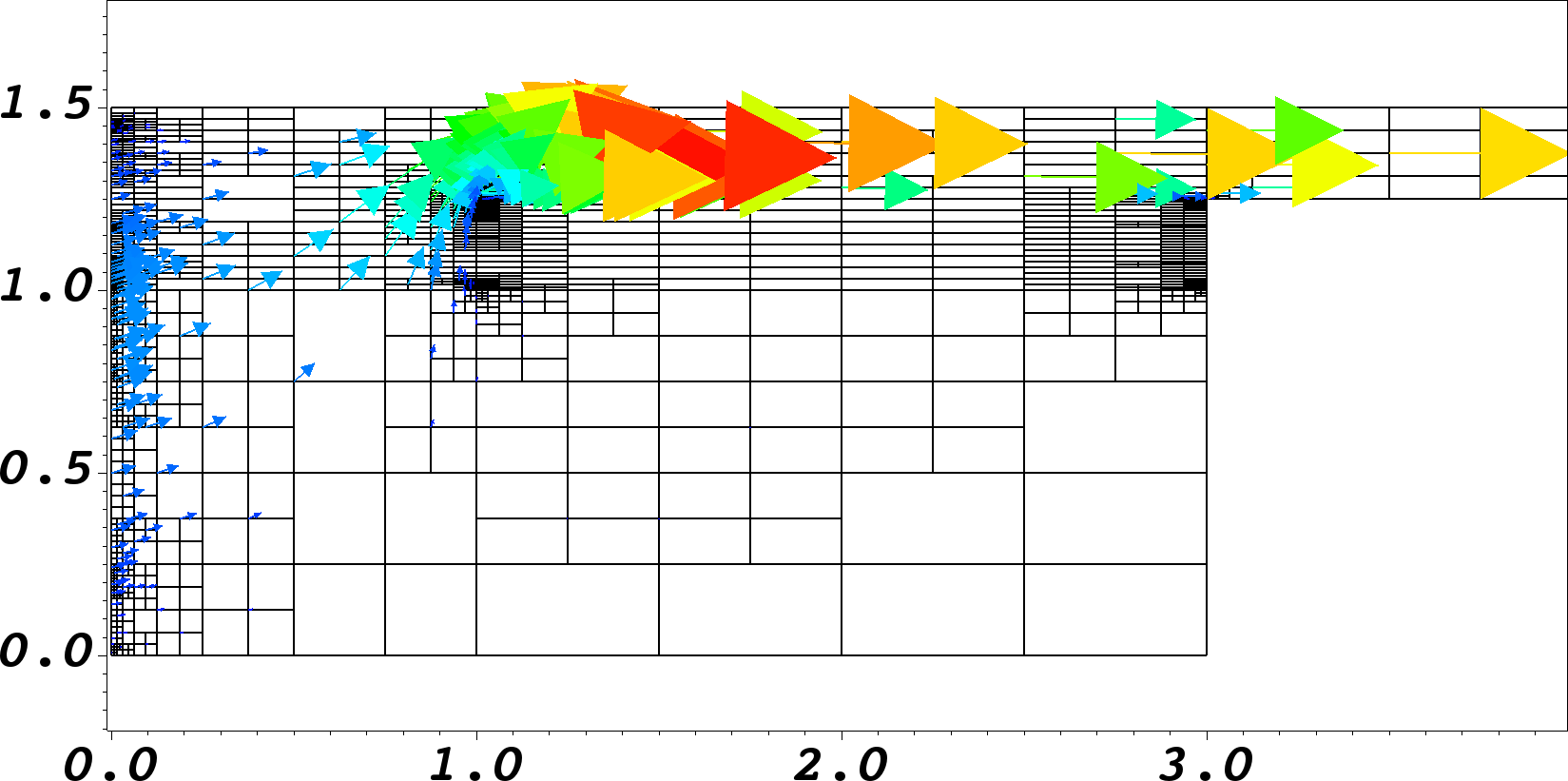}\qquad
  \includegraphics[width=78mm]{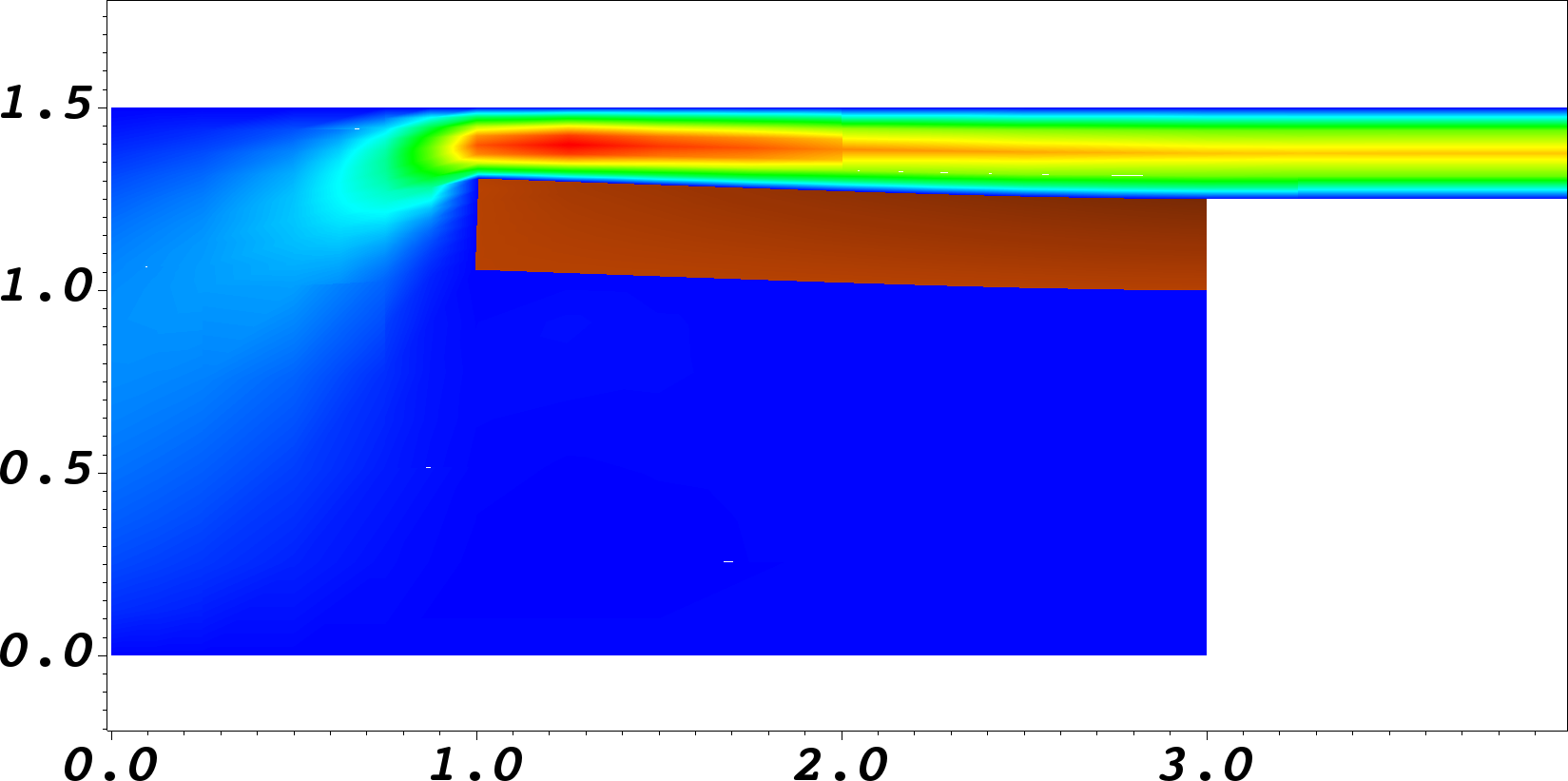}
  \caption{Example 2: locally refined mesh with vector plot
    of flow field (left) and flow field with elastic solid (brown)
    in the deformed configuration $\Omega$.}
  \label{pic_ex_3a}
\end{figure}

\paragraph{Configuration}
The geometrical data is sketched in \cite[Fig. 5.1, right subfigure]{RiWi13}
and also in the numerical results in Figure \ref{pic_ex_3a}.

\paragraph{Boundary conditions}
On the outer boundary $\partial\hat\Omega$ we work with homogeneous
Dirichlet conditions of the displacements. On the left inflow
boundary $\partial\hat\Omega\tsb{left}$ and right small channel
outflow $\partial\hat\Omega\tsb{right}$ pressure conditions are prescribed:
\begin{align*}
\hat g &= 0.2 I - \hrho_f \nu_f \hJ(\hF^{-T}\hnabla\hv_f^T\hn_f)\hF^{-T} \quad\text{on } \hGamma\tsb{left},\\
\hat g &= - \hrho_f \nu_f \hJ(\hF^{-T}\hnabla\hv_f^T\hn_f)\hF^{-T} \quad\text{on } \hGamma\tsb{right},
\end{align*}
where $I$ is the identity matrix in $\mathbb{R}^{2\times 2}$.
This means that we prescribe a pressure of \SI{0.2}{Pa} on the left
inflow part $\hGamma\tsb{left}$ and
zero pressure on the small channel outlet $\hGamma\tsb{right}$.
The second term is a correction
term due to the symmetric fluid stress tensor due to the
so-called do-nothing condition \cite{HeRaTu96}.
On the remaining outer boundaries, we prescribe $\hv=0$ (homogeneous Dirichlet conditions).

\paragraph{Parameters}

We use the fluid density $\varrho_f = \SI{1000}{kg m^{-3}}$ and kinematic
viscosity $\nu_f = \SI{0.001}{m^2 s^{-1}}$. The elastic
solid is characterized by
the density $\varrho_s = \SI{1000}{kg m^{-3}}$ and Poisson's ratio $\nu_s = 0.4$. Furthermore,
we use the Lam\'e coefficient $\mu_s = \SI{500}{kg m^{-1} s^{-2}}$.

\paragraph{Goals}
The combined functional consists
of a line integral evaluation $\Jdrag := F_D$ (see definition
below \eqref{drag_lift_forces}) and a point evaluation $J_2(\hp)=\hp(2.0,0.5)$
and reads
\[
J_c(\hU) = w_1 \Jdrag(\hU) + w_2 J_2(\hp).
\]
The weights are $w_1,w_2\in\mathbb{R}$ and $\omega_1=\omega_2=1$. We recall that $\hU := (\hv,\hu,\hp)$.
The reference value is computed on a sufficiently refined mesh:
\[
  J_c(\hat U) = \num[round-mode=places,round-precision=5]{2.7072783350606711e-01}.
\]

\begin{figure}[H]
  \centering
  \includegraphics[width=78mm]{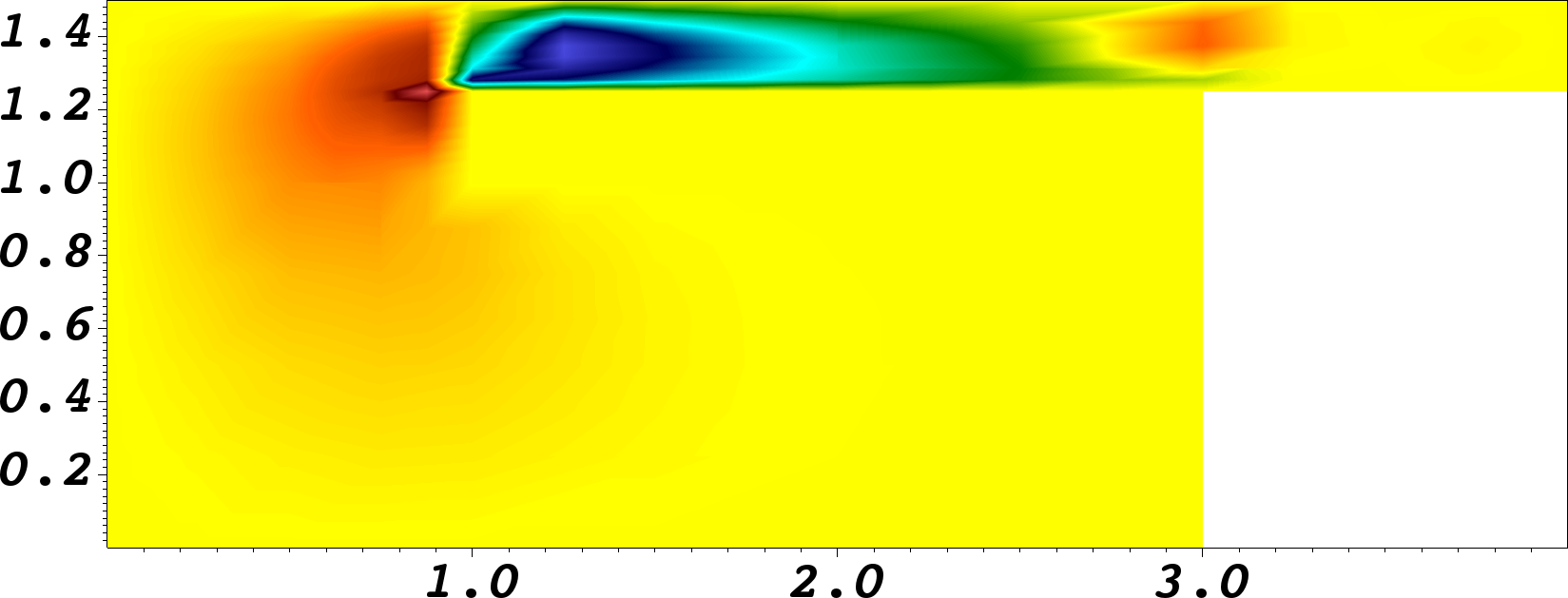}\qquad
  \includegraphics[width=78mm]{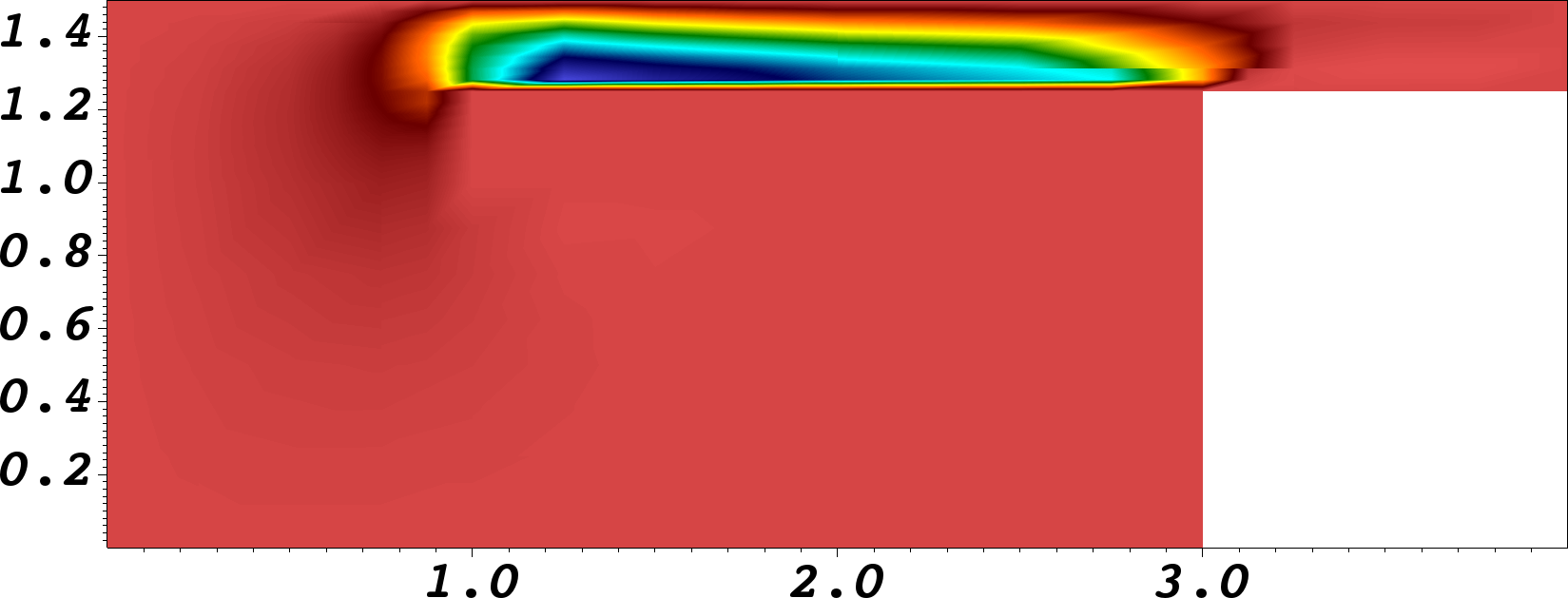}
  \caption{Example 2: adjoint displacement solutions $\hat z^{u_x}$
    and $\hat z^{u_y}$.}
  \label{pic_ex_3b}
\end{figure}

\paragraph{Results and discussion}
\ifcase0
Figures \ref{pic_ex_3a} and \ref{pic_ex_3b} display
the adaptively refined mesh, the elastic bar, the flow field
and the adjoint displacement solutions.
The adjoint solutions clearly indicate a strong influence
at the tip of the elastic bar and the small channel.
\or
Our graphical results are displayed in Figures \ref{pic_ex_3a} and
\ref{pic_ex_3b}. Therein, the adaptively refined mesh,
the elastic bar and the flow field are shown. Furthermore, the
adjoint displacement solutions are shown, which clearly indicate
a high influence at the tip of the elastic bar and the small channel.
\fi

The final values of the two goal functionals and the combined goal functional
are provided in Table~\ref{tab:ex2-final-values}.
\begin{table}[H]
  \centering
  \sisetup{round-mode=places,round-precision=3,table-format=1.3e+1}
  \begin{tabular}{lS}
    \toprule
    Drag     & 7.0815882385101531e-02 \\
    Pressure & 2.0000997093621237e-01 \\
    $J_c$    & 2.7082585332131393e-01 \\
    \bottomrule
  \end{tabular}
  \caption{Example 2: Final values of goal functionals.}
  \label{tab:ex2-final-values}
\end{table}

The true error, estimated error $\eta_h$ and corresponding
indices behave as given in Table~\ref{tab:ex2-levels}.
Again, we compare adaptive and uniform mesh refinement.
On refinement level $l=4$ we need approximately $10$ times more degrees of
freedom to achieve a comparable true error. However, uniform mesh refinement yields
better values in the estimator, which confirms our previous interpretation that adaptive
refinement is sensitive in terms of the error estimator.

\begin{table}[ht]
  \centering
  \sisetup{table-number-alignment=center,table-format=1.2e+1}
  \begin{tabular}{rrSSS
    %S[scientific-notation=fixed,fixed-exponent=0,table-format=3.3]
    %S[scientific-notation=fixed,fixed-exponent=0,table-format=3.2]
    SS}
    \toprule
    Level $l$ & Dofs&{$|J_c(\hat U)-J_c(\hat U_h)|$}&{$|\eta_h|$}&{$\sum_i|\eta_i|$}&$\Ieff$&$\Iind$ \\
    \midrule
    1 & 2125  & 8.48e-04 & 5.59e-04 & 9.57e-04 & 6.59e-01 & 1.13e+00 \\
    2 & 5508  & 3.49e-04 & 8.23e-05 & 5.79e-04 & 2.36e-01 & 1.66e+00 \\
    3 & 10674 & 1.39e-04 & 1.74e-03 & 2.09e-03 & 1.25e+01 & 1.50e+01 \\
    4 & 16383 & 4.62e-05 & 1.15e-03 & 1.36e-03 & 2.49e+01 & 2.94e+01 \\
    5 & 25994 & 2.89e-06 & 2.16e-03 & 2.30e-03 & 7.47e+02 & 7.94e+02 \\
    6 & 36596 & 9.80e-05 & 1.34e-03 & 1.47e-03 & 1.37e+01 & 1.50e+01 \\ \midrule
    \multicolumn7c{Uniform mesh refinement} \\
    \midrule
1 & 2125  &  8.48e-04   &     5.59e-04  &      9.57e-04   &     6.59e-01  &      1.13e+00\\
2 & 8053  &  3.50e-04   &     1.95e-04  &      4.29e-04   &     5.57e-01  &      1.22e+00\\
3 & 31333 &  1.36e-04   &     6.38e-05  &      2.16e-04   &     4.69e-01  &      1.59e+00\\
4 & 123589&  4.14e-05   &     1.71e-05  &      1.16e-04   &     4.12e-01  &      2.80e+00\\
    \bottomrule
  \end{tabular}
  \caption{Example 2: Degrees of freedom, true error, estimator and indices.}
  \label{tab:ex2-levels}
\end{table}
We observe sufficient decrease in the true error and estimated error as well
as the indicator index. However, there is a large difference of about one order
of magnitude resulting in an overestimation of both the effectivity
and indicator indices. The reason for this is the specific
setting, which is very sensitive to the ALE transformation and also to the
relatively large influence of the adjoint. %
%For the latter, we would expect better result
%by working with the adjoint error part, too.

%%%%%%%%%%%%%%%%%%%%%%%%%%%%%%%%%%%%%%%%%%%%%%%%%%%%%%%%%%%%%%%
\subsection{Example 3: FSI-1 benchmark}
This configuration is taken from \cite{HrTu06b}. Our own results using
uniform mesh refinement are provided in \cite{Wi13_fsi_with_deal,RiWi10}.
The extension to multiple goal functionals is novel.

% Internal comments
%Hron/Turek-benchmark paper: siehe
%p. 12 Table 13.
%
%Fuer den FSI-1 gibt es insgesamt 4 goal functionals, darunter zwei Mal Punktauswertungen der
%Flaggenspitze.
%
%Viel mehr faellt mir auch nicht ein und sollte auch reichen.
%
%Man koennte drei Settings rechnen:
%
%1. Nur single goal, z.B. lift
%
%2. Zwei goal functionals: drag und $u_y$
%
%3. Alle vier goal functionals: drag, lift, ux, uy
%
%Im Vergleich von 2. und 3. wuerde ich beim adaptiven Gitter kaum Unterschiede erwarten, da
%die goal functionals sich jeweils sehr ähneln. Aber das waere interessant zu sehen.

\paragraph{Configuration}

The computational domain has length $L=\SI{2.5}{m}$ and height $H=\SI{0.41}{m}$. The circle center
is positioned at $C=(\SI{0.2}{m},\SI{0.2}{m})$ with radius $r=\SI{0.05}{m}$. The elastic beam has length
$l=\SI{0.35}{m}$ and height $h=\SI{0.02}{m}$. The right lower end is positioned at
$(\SI{0.6}{m},\SI{0.19}{m})$, and
the left end is attached to the circle.

\begin{figure}[h]
\centering
{\includegraphics[width=10cm]{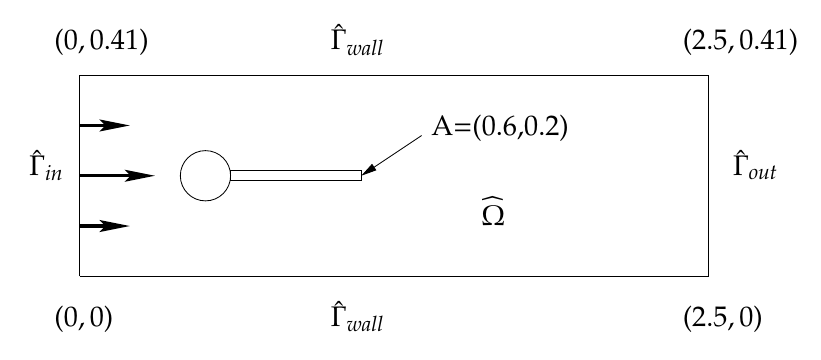}}
\caption{Flow around cylinder with elastic beam with
circle-center $C=(0.2,0.2)$ and radius $r=0.05$.}
\label{configuration_csm_and_fsi_2D}
\end{figure}

Control points $A(t)$ (with $A(0) = (0.6,0.2)$) are fixed at the
trailing edge of the structure, measuring $x$- and $y$-deflections of the beam.

\paragraph{Boundary conditions}

A parabolic inflow velocity profile is given on $\hat\Gamma\tsb{in}$ by
\begin{align*}
  v_f (0,y) = 1.5\bar{U} \frac{4y(H-y)}{H^2} , \quad \bar{U} = \SI{0.2}{ms^{-1}}.
\end{align*}
On the outlet $\hat\Gamma\tsb{out}$ the `do-nothing' outflow condition \cite{HeRaTu96}
is imposed
which leads to zero mean value of the pressure at this part of the boundary.
The displacements are fixed (homogeneous Dirichlet conditions) around the outer
boundary and the cylinder. On $\hGamma\tsb{wall}$ no-slip conditions for flow
(homogeneous Dirichlet for velocities $\hv$) are prescribed.

\paragraph{Parameters}

We use the fluid density $\varrho_f = \SI{1000}{kg m^{-3}}$ and kinematic
viscosity $\nu_f = \SI{0.001}{m^2 s^{-1}}$. The elastic
solid is characterized by
the density $\varrho_s = \SI{100}{kg m^{-3}}$ and Poisson's ratio $\nu_s = 0.4$. Furthermore,
we use the Lam\'e coefficient $\mu_s = \SI{0.5e6}{kg m^{-1} s^{-2}}$.

\paragraph{Quantities of interest of the original benchmark problem}
In the benchmark configuration \cite{HrTu06b} four quantities of interest
were evaluated:

\begin{enumerate}
\item[1)] $x$- and $y$-deflection of the beam at $A(t)$.
\item[2)] The forces exerted by the fluid on the whole body,
i.e., drag force $F_D$ and lift force $F_L$
on the rigid cylinder and the elastic beam. They form a closed path in which
the forces can be computed with the help of line integration.
The formula is evaluated on the fixed reference domain $\hat\Omega$ and reads:
\begin{align}\label{drag_lift_forces}
(F_D , F_L) &= \int_{\hat S} \hat J\hat\sigma\tsb{all}\hat F^{-T} \cdot \hat n\, d\hat s
= \int_{\hat S(\text{circle})} \hat J \hat\sigma_f \hat F^{-T} \cdot \hat n_f \, d\hat s +
\int_{\hat S(\text{beam})} \hat J\hat\sigma_f \hat F^{-T} \cdot \hat n_f \, d\hat s.
\end{align}
\end{enumerate}

\paragraph{Goals}
We consider three goal functionals simultaneously.
The previous benchmark quantities of interest all yield local refinements
around the FSI interface, when goal-oriented error control is employed.
In order to highlight
more clearly that three really distinct functionals can be controlled, we choose
(1) the above drag functional, (2) away from the FSI interface, a pressure
point evaluation in $J_2(\hp) := \hp(1.5,0.3)$, and (3) the flux evaluation on the boundary $\hGamma\tsb{out}$, i.e.,
$\int_{\hGamma\tsb{out}} \hv_f\cdot \hat n\, ds$.
Thus, the combined functional reads
\[
J_c(\hU) = w_1 \Jdrag(\hU) + w_2 J_2(\hp) + w_3 J_3(\hv).
\]
The weights are $w_1,w_2,w_3\in\mathbb{R}$, $\omega_1=\omega_2=\omega_3=1$,
and we recall that $\hU := (\hv,\hu,\hp)$.
The reference value is computed on a sufficiently refined mesh:
\[
J_c(\hat U) = \num[round-mode=places,round-precision=5]{3.1205264275814216e+01}.
\]

\begin{figure}[ht]
  \centering
  \includegraphics[width=14cm]{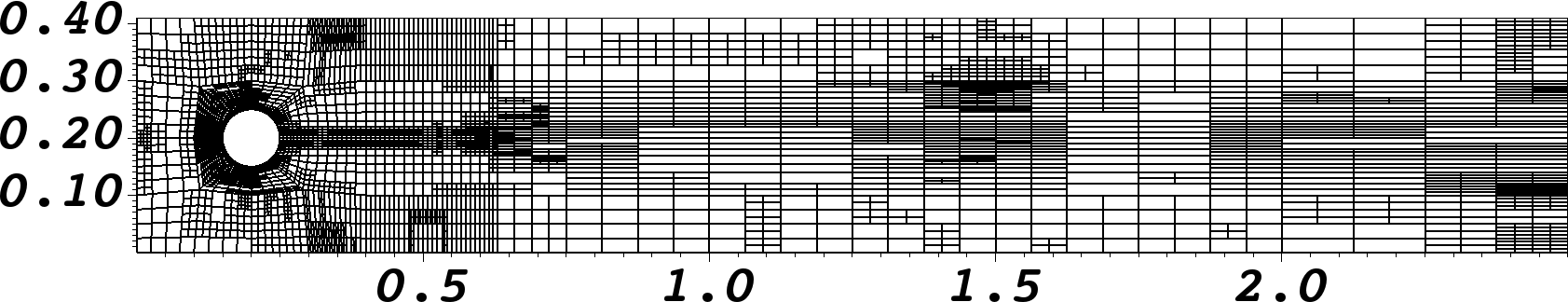}\\[1ex]
  \includegraphics[width=14cm]{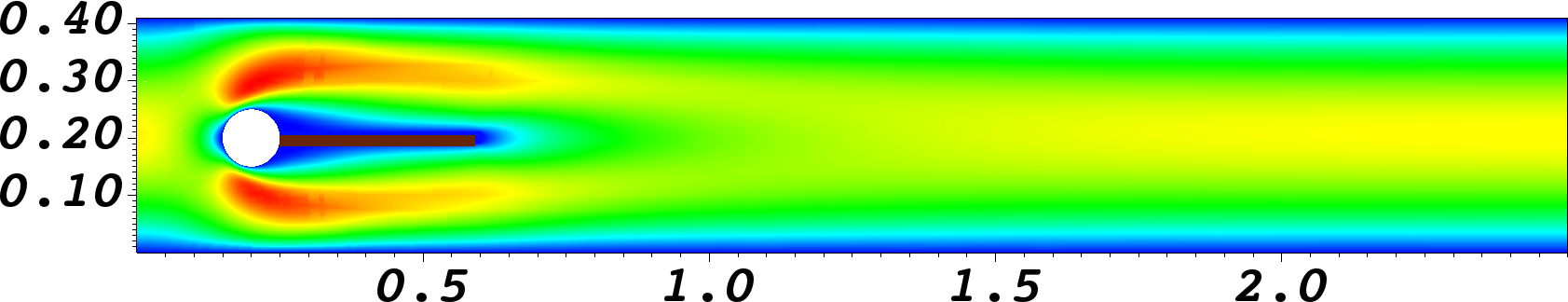}\\[1ex]
  \includegraphics[width=14cm]{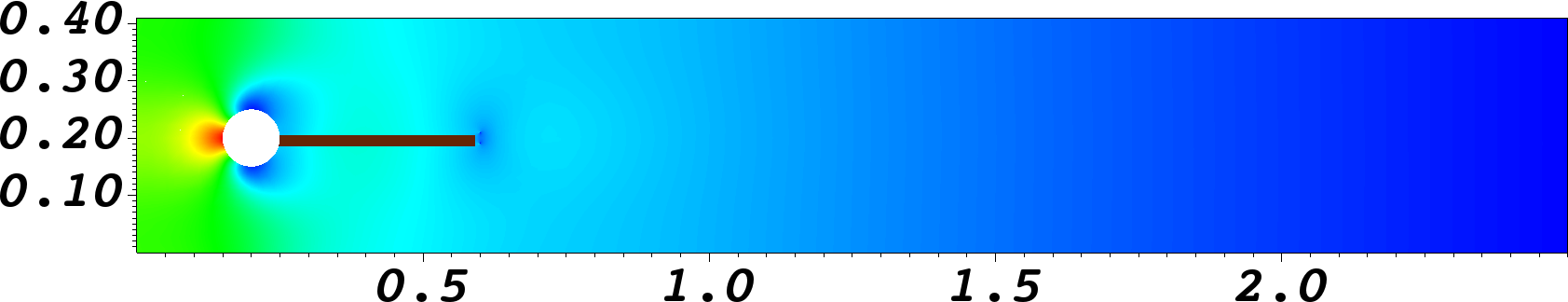}
  \caption{From top to bottom: adaptively refined mesh,
    $\hv_s$ velocity field, and the pressure field
    after three adaptive refinement steps.
    The elastic beam is colored in blue in the middle figure.
    In the top figure, we observe adaptive mesh refinement
    mainly around the interface between the elastic beam
    and the surrounding fluid and also in the
    pressure point evaluation at $J_2(\hp) := \hp(1.5,0.3)$.}
  \label{pic_ex_1a}
\end{figure}

\begin{figure}[H]
  \centering
  \includegraphics[width=14cm]{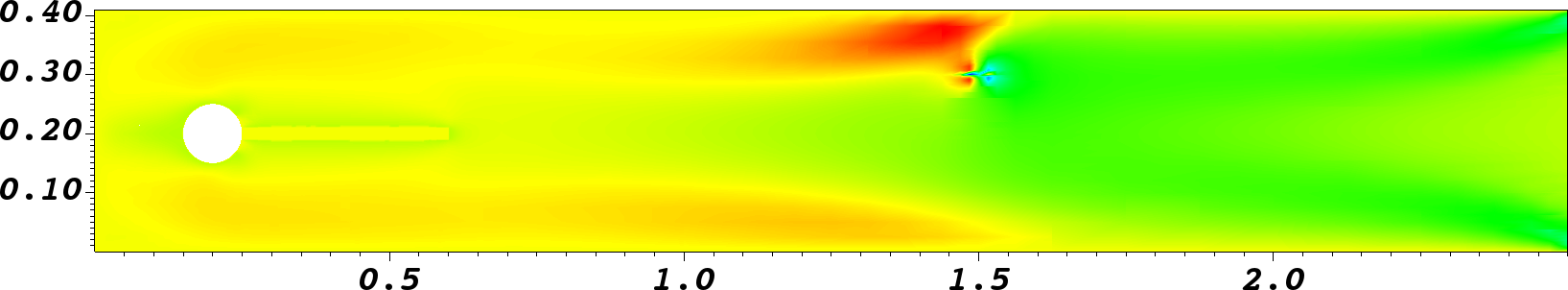}\\[1ex]
  \includegraphics[width=14cm]{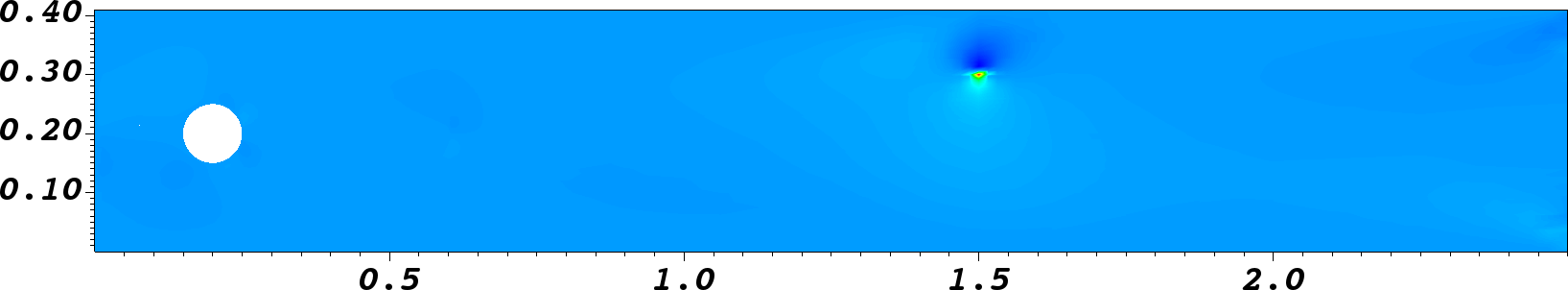}\\[1ex]
  \includegraphics[width=14cm]{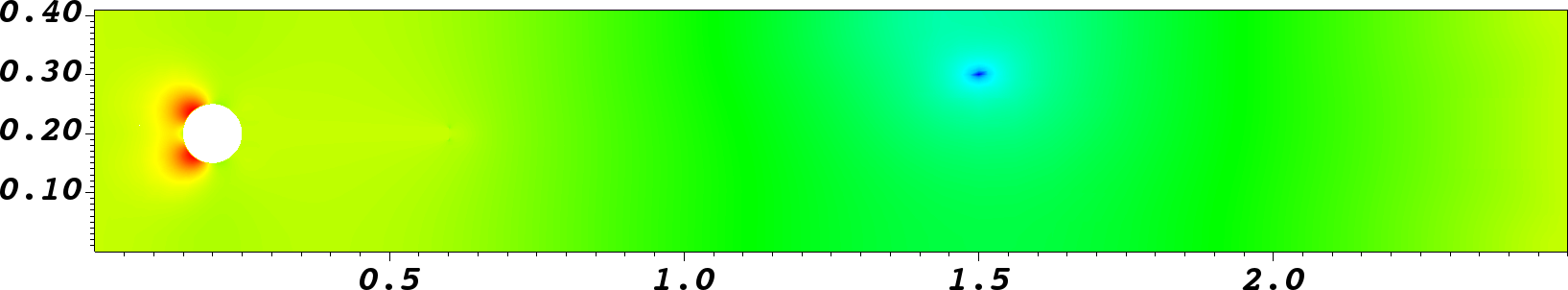}\\[1ex]
  \includegraphics[width=14cm]{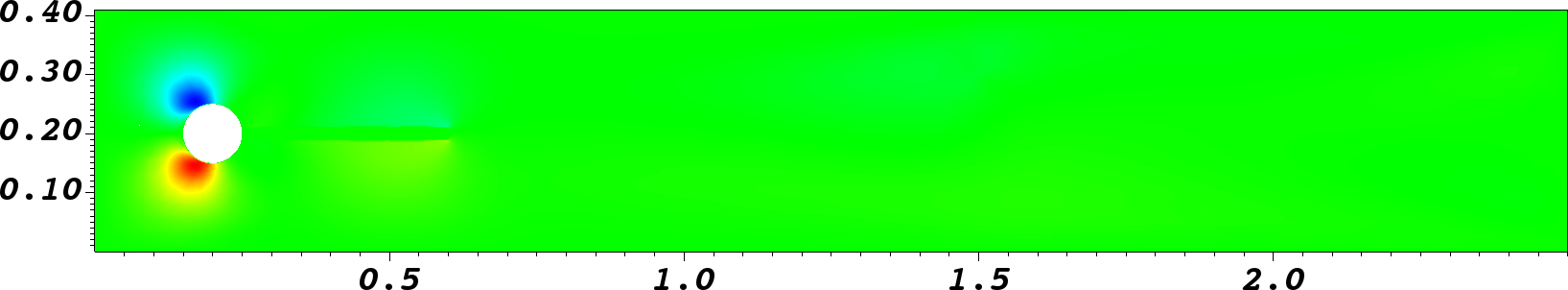}\\[1ex]
  \includegraphics[width=14cm]{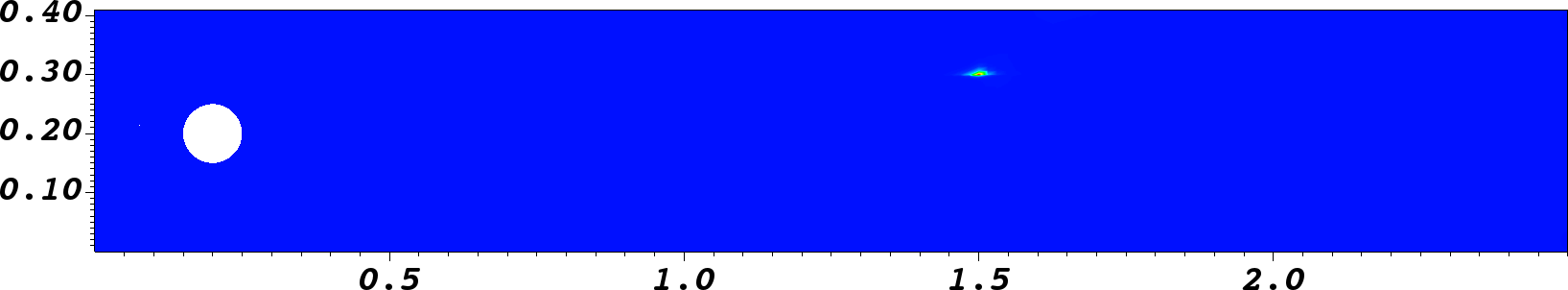}
  \caption{From top to bottom: adjoint solutions for
    $\hz^{\hv_x},\hz^{\hv_y},\hz^{\hu_x},\hz^{\hu_y}$ and $\hz^{\hp}$.}
  \label{pic_ex_1b}
\end{figure}

\paragraph{Results and discussion}
In Figure \ref{pic_ex_1a}, the final adaptive
mesh and primal solution fields are shown. All adjoint solution fields
are shown in Figure \ref{pic_ex_1b}.

Independently of the specific multi-goal framework, we first notice
that the obtained FSI benchmark values agree well with the literature values
\cite{HrTu06b,Wi13_fsi_with_deal}, see Table~\ref{tab:ex3-goal-final-a}.
\begin{table}[ht]
  \centering
  \sisetup{round-mode=places,round-precision=3,table-format=1.3e+1}
  \begin{tabular}{lS}
    \toprule
    DisX & 2.2657479709296053e-05 \\
    DisY & 8.2001891962646791e-04 \\
    Drag & 1.5351737833128903e+01 \\
    Lift & 7.3885947240664507e-01 \\
    \bottomrule
  \end{tabular}
  \caption{Example 3: Final values of $x$- and $y$-deflection
    of the beam at $A(t)$, drag, and lift.}
  \label{tab:ex3-goal-final-a}
\end{table}

The final drag, pressure value and flux values (multiple goal functionals)
and the
combined goal functional are given in Table~\ref{tab:ex3-goal-final-b}.
\begin{table}[H]
  \centering
  \sisetup{round-mode=places,round-precision=3,table-format=1.3e+1}
  \begin{tabular}{lS}
    \toprule
    Drag     & 1.5351737833128903e+01 \\
    Pressure & 1.5766176006021523e+01 \\
    Flux     & 8.1999999975009522e-02 \\
    $J_c$    & 3.1199913839125436e+01 \\
    \bottomrule
  \end{tabular}
  \caption{Example 3: Final values of the multigoal functionals.}
  \label{tab:ex3-goal-final-b}
\end{table}

With regard to our three goal functionals and error estimator,
we obtained the results given in Table~\ref{tab:ex3-levels}.
Therein, at the second and third refinement levels we need
two and four times more degrees of freedom, respectively, to achieve comparable
true errors and estimators, which shows that our adaptive multigoal scheme
can achieve a significant reduction of computational cost.
\begin{table}[H]
  \centering
  \sisetup{table-number-alignment=center,table-format=1.2e+1}
  \begin{tabular}{rrSSS
    S[scientific-notation=fixed,fixed-exponent=0,table-format=1.3]
    S[scientific-notation=fixed,fixed-exponent=0,table-format=1.2]}
    \toprule
    Level $l$ & Dofs&{$|J_c(\hat U)-J_c(\hat U_h)|$}&{$|\eta_h|$}&{$\sum_i|\eta_i|$}&$\Ieff$&$\Iind$ \\
    \midrule
    1 & 13310  & 2.73e-01 & 1.42e-01 & 5.41e-01 & 5.20e-01 & 1.98e+00 \\
    2 & 27193  & 7.87e-02 & 4.46e-02 & 1.69e-01 & 5.66e-01 & 2.15e+00 \\
    3 & 54893  & 2.05e-02 & 1.24e-02 & 6.11e-02 & 6.07e-01 & 2.98e+00 \\
    4 & 109909 & 5.35e-03 & 4.42e-03 & 2.88e-02 & 8.26e-01 & 5.39e+00 \\ \midrule
    \multicolumn7c{Uniform mesh refinement} \\
    \midrule
1 & 13310  & 2.73e-01  &      1.42e-01  &      5.41e-01  &      5.20e-01 &       1.98e+00\\
2 & 52052  & 7.48e-02  &      4.52e-02  &      1.60e-01  &      6.05e-01 &       2.14e+00\\
3 & 205832 & 2.04e-02  &      1.32e-02  &      5.64e-02  &      6.49e-01 &       2.76e+00\\
    \bottomrule
  \end{tabular}
  \caption{Example 3: Degrees of freedom, true error, estimator and indices.}
  \label{tab:ex3-levels}
\end{table}
Therein, we observe well that the estimated error decreases by two orders of magnitude
and $\Ieff$ and $\Iind$ show a nice behavior. In view of the complexity
of the problem statement dealing with a nonlinear, coupled fluid-structure
interaction system and only working with the primal error estimator, these are excellent results
which show that the multigoal-technology performs well for this example.
%Moreover, both functional values converge to limit values.
%Overall, the multigoal-technology performs well for this example.

%%%%%%%%%%%%%%%%%%%%%%%%%%%%%%%%%%%%%%%%%%%%%%%%%%%%%%%%%%%
\section{Conclusions}
In this work, we developed multigoal-oriented a posteriori
error control for stationary fluid-structure interaction.
Specifically, the simultaneous control of several quantities of interest
for such multiphysics problems may be
required in practical applications. The focus in this work was on
prototype settings for verification of our proposed multigoal framework.
Three numerical tests, namely an elastic lid-driven
cavity, a chamber with elastic solid and the FSI-1 benchmark
were adopted to study the performance of our methodology.
Therein, the lid-driven
cavity (Example 1) and the FSI-1 benchmark (Example~3) yield
good error reductions, good estimators and therefore
good effectivity indices. The resulting adaptive meshes localize
well the different goal functionals. Moreover, different weights 
were investigated in Example 1. In Example 2, the error reductions
are also good, but the effectivity index shows an overestimation of
about one order magnitude. This test is somewhat challenging due to the
large solid displacement and the influence of the adjoint solution.
Two (Examples 1 and 2) and three (Example 3) goal functionals
were studied and the overall performance is excellent in view of the
complexity of the governing fluid-structure interaction system.
A future (challenging) extension is the development of a
multigoal framework for time-dependent fluid-structure interaction
problems.

%%%%%%%%%%%%%%%%%%%%%%%%%%%%%%%%%%%%%%%%%%%%%%%%%%%%%%%%%%%
\section*{Acknowledgments}
This work has been supported by the DAAD in the project
`A new passage to India' between the Leibniz University Hannover
and IIT Indore. Furthermore, the second and last author are affiliated
to the Cluster of Excellence PhoenixD (EXC 2122, Project ID 390833453).
Moreover, the second author acknowledges support from
the Austrian Science Fund (FWF) under the grant
P-29181
`Goal-Oriented Error Control for Phase-Field Fracture Coupled
to Multiphysics Problems' at the beginning of this work in Linz.
The third author is funded by the Deutsche Forschungsgemeinschaft
(DFG, German Research Foundation) -- SFB1463 -- 434502799.

%%%%%%%%%%%%%%%%%%%%%%%%%%%%%%%%%%%%%%%%%%%%%%%%%%%%%%%%%%%

%\bibliographystyle{abbrv}
%\bibliography{lit}

\end{document}